\documentclass{amsart}
\usepackage[a4paper,margin=2.5cm]{geometry}
\usepackage{amssymb, colordvi}
\usepackage{multirow}
\usepackage[pdftex]{graphicx,color}
\usepackage{xfrac, faktor}
\usepackage{mathtools}

\numberwithin{equation}{section}
\newtheorem{thm}{Theorem}
\numberwithin{thm}{section}
\newtheorem{prop}[thm]{Proposition}
\newtheorem{lemma}[thm]{Lemma}

\newtheorem{cor}[thm]{Corollary}
\theoremstyle{definition}

\newtheorem{example}[thm]{Example}
\newtheorem{remark}[thm]{Remark}

\newcommand{\sgnvar}{{\operatorname{sgnvar}}}
\newcommand{\sgn}{{\operatorname{sgn}}}
\newenvironment{ex}{\begin{example}\rm}{\end{example}}
\newenvironment{rem}{\begin{remark}\rm}{\end{remark}}
\newcounter{FNC}[page]
\def\fauxfootnote#1{{\addtocounter{FNC}{2}$^\fnsymbol{FNC}$%
     \let\thefootnote\relax\footnotetext{$^\fnsymbol{FNC}$#1}}}

\newcommand{\N}{\mathbb{N}}

\newcommand{\R}{\mathbb{R}}
\newcommand{\Z}{\mathbb{Z}}

\newcommand{\A}{{\mathcal A}}

\newcommand{\Ci}{{\mathcal C}^{\infty}}

\title[Budan--Fourier theorem for circuit polynomial systems]{Budan--Fourier theorem for circuit polynomial systems}

\author{Fr\'ed\'eric Bihan}
\address{Laboratoire de Math\'ematiques\\
         Universit\'e Savoie Mont Blanc\\
         73376 Le Bourget-du-Lac Cedex\\
         France}

\email{Frederic.Bihan@univ-savoie.fr}
\urladdr{http://www.lama.univ-savoie.fr/~bihan/}


\begin{document}

\begin{abstract}
We prove a Budan--Fourier theorem for linear combinations of analytic functions whose successive Wronskians are positive: the number of zeros of $g=\sum_{k=1}^s c_kh_k$ in an interval $(a,b]$, counted with multiplicities, is bounded by the drop of the number of sign changes of the sequence of Wronskians $W(h_1,\ldots,h_{k-1},g)$ between $a$ and $b$, and the difference is even. This refines the bound $s-1$ of Koiran, Portier and Tavenas, and we show that it applies, after a suitable renumbering, to any sequence of functions satisfying Descartes' rule of signs.
We then study the functions arising from polynomial systems supported on circuits, that is, systems of $n$ equations in $n$ variables with $n+2$ monomials. For these functions, we show that the Budan--Fourier bound never exceeds Descartes' bound, and we compute it at the endpoints of the interval of definition. As an application, we obtain a refinement of the optimal Descartes' rule of signs for circuits of Bihan, Dickenstein and Forsg{\aa}rd: we bound the number of positive solutions whose logarithm lies in explicit polytopes, these bounds add up to at most the bound of Bihan, Dickenstein and Forsg{\aa}rd, and they determine the parity of the number of positive solutions.
\end{abstract}

\maketitle
\section{Introduction}

Descartes' rule of signs asserts that the number of positive roots of a real univariate polynomial $f=c_0+c_1x+\cdots+c_dx^d$, counted with multiplicities, is at most the number $\sgnvar(c_0,\ldots,c_d)$ of sign changes in the sequence of its coefficients, and that the difference is even. The Budan--Fourier theorem refines it: for $a<b$, the number of roots of $f$ in $(a,b]$ is at most $\sgnvar(f(a),f'(a),\ldots,f^{(d)}(a))-\sgnvar(f(b),f'(b),\ldots,f^{(d)}(b))$, and the difference is even. Both results have been extended far beyond polynomials. A sequence of analytic functions $h_1,\ldots,h_s$ on an interval $I$ is said to satisfy Descartes' rule of signs if the number of zeros of $\sum_k c_kh_k$ in $I$ never exceeds $\sgnvar(c_1,\ldots,c_s)$; P\'olya and Szeg\H{o} characterized such sequences in terms of the signs of their Wronskians \cite{P-S}. Coste, Lajous-Loaeza, Lombardi and Roy proved a Budan--Fourier theorem in which the derivatives are replaced by the operators $g \mapsto (g/f_i)'$ associated with a sequence of positive functions $f_i$ \cite{C}. Koiran, Portier and Tavenas used Wronskians to bound the number of real zeros of a linear combination $g=\sum_{k=1}^s c_kh_k$: if none of the Wronskians $W(h_1,\ldots,h_k)$, $k=1,\ldots,s$, vanishes on $I$, then $g$ has at most $s-1$ zeros in $I$, counted with multiplicities, and in general the number of zeros of $g$ is controlled by the number of zeros of these Wronskians \cite{KT}.

\subsection*{A Budan--Fourier theorem with Wronskians}
Our first result combines these approaches. Let $h_1,\ldots,h_s$ be analytic functions on an open interval $I$ such that $W(h_1,\ldots,h_k)>0$ on $I$ for $k=1,\ldots,s$, let $g=\sum_{k=1}^s c_kh_k$ with $c_s \neq 0$, and consider the sequence of functions
$$H(g)=\big(H_1(g),\ldots,H_s(g)\big), \qquad H_k(g)=W(h_1,\ldots,h_{k-1},g).$$
Theorem \ref{application} asserts that for $a<b$ in $I$,
$$n(g;(a,b]) \leq \sgnvar H(g)(a)-\sgnvar H(g)(b),$$
where $n(g;(a,b])$ is the number of zeros of $g$ in $(a,b]$ counted with multiplicities, and that the difference is even. Since $H(g)$ has $s$ terms, this refines the bound $s-1$ of \cite{KT}. For $h_k=x^{k-1}$ we have $H_k(g)=c\,g^{(k-1)}$ with $c>0$, and we recover the classical Budan--Fourier theorem; for $h_k=x^{a_k}$ with real exponents $a_1<\cdots<a_s$ we recover Descartes' rule of signs on $(0,+\infty)$ (Example \ref{ex:monomialsW}). The proof relies on the Frobenius--P\'olya factorization of the operator defined by $h_1,\ldots,h_s$ and on the generalized Budan--Fourier theorem of \cite{C}. We also show that any sequence of analytic functions satisfying Descartes' rule of signs can be renumbered, and multiplied by a common sign, so that all the Wronskians $W(h_1,\ldots,h_k)$ become positive (Proposition \ref{prop:renumbering}); in general a non-trivial permutation is needed. Consequently, Theorem \ref{application} applies to all the families of functions satisfying Descartes' rule of signs given in \cite[Part V]{P-S}, for instance exponential functions $e^{\lambda_1x},\ldots,e^{\lambda_sx}$ with $\lambda_1<\cdots<\lambda_s$ (Laguerre's rule of signs for exponential sums), monomials $x^{a_1},\ldots,x^{a_s}$ with real exponents on $(0,+\infty)$, their confluent versions such as $x^je^{\lambda x}$ and $x^a(\log x)^j$, or functions $K(x,t_1),\ldots,K(x,t_s)$ defined by a totally positive kernel $K$ \cite{Karlin}. Apart from monomials (Example \ref{ex:monomialsW}), we do not work out these families in this paper; see Remark \ref{rem:powers} for a few comments.

\subsection*{A refinement of Descartes' rule of signs}
The initial motivation of this paper was to refine the optimal Descartes' rule of signs for polynomial systems supported on circuits of Bihan, Dickenstein and Forsg{\aa}rd \cite{BDF}, recalled below. By Gale duality, this amounts to refining Descartes' rule of signs for a family of functions built from linear functions, and our main application of Theorem \ref{application} is such a refinement. Let $p_j(y)=a_j+b_jy$, $j=1,\ldots,s$, be linear functions with $a_ib_j-a_jb_i<0$ for $i<j$, and assume that the interval $I$ where all the $p_j$ are positive is non-empty and bounded, say $I=(\alpha,\beta)$ (the unbounded case reduces to this one, see Remark \ref{rem:unbounded}). Then $\alpha$ is the root of $p_1$ and $\beta$ is the root of $p_s$.

\smallskip

\noindent\emph{Functions arising from circuits.} Consider the functions
$$h_i=\frac{b_i}{p_i}-\frac{b_{i+1}}{p_{i+1}}, \qquad i=1,\ldots,s-1,$$
which appear in the study of polynomial systems supported on circuits \cite{BD,BDF} and satisfy Descartes' rule of signs on $I$: for $g=\sum_{i=1}^{s-1}c_ih_i$, the number $n(g;I)$ of roots of $g$ in $I$ is at most $\sgnvar(c_1,\ldots,c_{s-1})$. Let $H(g)=(H_1(g),\ldots,H_{s-1}(g))$, with $H_k(g)=W(h_1,\ldots,h_{k-1},g)$, be the Budan--Fourier sign sequence. We show that it never exceeds Descartes' bound, $\sgnvar H(g)(y) \leq \sgnvar(c_1,\ldots,c_{s-1})$ for every $y \in I$ (Proposition \ref{prop:BFleqD}); the proof uses an induction on $s$ and the invariance of $H(g)$ under projective changes of variable. Computing the limits of $\sgnvar H(g)$ at the endpoints of $I$, we obtain a number $V_\alpha(c)$, computed from the values of the functions $h_i$ at $\alpha$, such that
$$n(g;I) \;\leq\; V_\alpha(c) \;\leq\; \sgnvar(c_1,\ldots,c_{s-1}), \qquad n(g;I) \equiv V_\alpha(c) \equiv \sgnvar(c_1,\ldots,c_{s-1}) \mod 2$$
(Corollary \ref{cor:VleqD}). The second inequality can be strict, so that we get a genuine refinement of Descartes' rule of signs, and the bound $V_\alpha(c)$ is not determined by the signs of the coefficients. For instance, for $s=4$ we get $V_\alpha(c)=\sgnvar\big(c_1,\ c_2h_2(\alpha)+c_3h_3(\alpha),\ c_3\big)$: if $c_1,c_3>0>c_2$ and $c_3h_3(\alpha)>|c_2|\,h_2(\alpha)$, then $g$ has no root in $I$, while Descartes' rule of signs allows two roots (Remark \ref{rem:comparisonDescartes}).

\smallskip

\noindent\emph{Negative powers of linear functions.} The same phenomenon occurs for negative powers of the linear functions $p_j$, and in that case all the quantities involved are explicit. Let $\gamma<0$ and $F=\sum_{j=1}^s c_jp_j^{\gamma}$ with $c_1c_s \neq 0$. The functions $p_1^\gamma,\ldots,p_s^\gamma$ satisfy Descartes' rule of signs on $I$, and we prove (Theorem \ref{thm:powers}) that
$$n(F;I) \;\leq\; V^\gamma_\alpha(c) \;\leq\; \sgnvar(c_1,\ldots,c_s), \qquad n(F;I) \equiv V^\gamma_\alpha(c) \equiv \sgnvar(c_1,\ldots,c_s) \mod 2,$$
where $V^\gamma_\alpha(c)=\sgnvar\big(c_1,\sigma_2(\alpha),\ldots,\sigma_s(\alpha)\big)$ is given by the explicit numbers
$$\sigma_k(\alpha)=\sum_{i=k}^s c_i \prod_{j=1}^{k-1}|a_jb_i-a_ib_j| \; p_i(\alpha)^{\gamma-k+1}, \qquad k=2,\ldots,s,$$
assumed to be non-zero. Again the second inequality can be strict: for $s=3$ we get $V^\gamma_\alpha(c)=\sgnvar\big(c_1,\ c_2p_2(\alpha)^\gamma+c_3p_3(\alpha)^\gamma,\ c_3\big)$, so that if $c_1,c_3>0>c_2$ and $c_3p_3(\alpha)^\gamma>|c_2|\,p_2(\alpha)^\gamma$, then $F$ has no root in $I$, while Descartes' rule of signs allows two roots. More generally, the Budan--Fourier sign sequence of $F$ at any point $y \in I$ is the sequence of signs of the explicit functions $\sigma_k(y)$, which does not involve Wronskians, and its number of sign changes never exceeds $\sgnvar(c_1,\ldots,c_s)$. This refinement fails in general for positive exponents (Remark \ref{rem:powers}).

\subsection*{Positive solutions of systems supported on circuits}
Let $\A=\{a_0,\ldots,a_{n+1}\} \subset \Z^n$ be a circuit and $C \in \R^{n \times (n+2)}$ a coefficient matrix of rank $n$, and denote by $n_\A(C)$ the number of positive solutions of the system $\sum_j c_{i,j}x^{a_j}=0$, $i=1,\ldots,n$. Bihan, Dickenstein and Forsg{\aa}rd proved the optimal bound $n_\A(C) \leq 1+\sgnvar(\mu)$, where $\mu$ is a sequence of partial sums of the coefficients of the affine relation of $\A$, ordered according to the signs of the maximal minors of $C$ \cite{BDF}, refining a previous bound of \cite{BD}. Their proof reduces, by Gale duality \cite{BS07}, to bounding the number of zeros of a function $G=\log g$ on an interval $\Delta_P$, and applies Descartes' rule of signs to $G'$, which is a linear combination of functions as above, together with Rolle's theorem.

Applying our Budan--Fourier theorem directly to $G$, we obtain Theorem \ref{thm:localized}. It provides an explicit sign sequence $S(y)$ such that, for any subinterval $(a,b] \subset \Delta_P$, the number of positive solutions $x$ such that $\log x$ lies in an explicit polyhedron $\mathcal{P}_{(a,b]} \subset \R^n$, defined by linear inequalities involving the differences of the exponents $a_i-a_j$ (a polytope when $[a,b] \subset \Delta_P$), is at most $S(a)-S(b)$, with an even difference. Moreover $S \leq 1+\sgnvar(\mu)$, so that the bounds attached to a subdivision of $\Delta_P$ add up to at most the bound of \cite{BDF}, and, under the hypotheses of \cite{BDF} ensuring parity, the limit of $S$ at the left endpoint of $\Delta_P$ gives a bound $1+V_\alpha \leq 1+\sgnvar(\mu)$ for $n_\A(C)$ which also determines its parity. In particular, for any compact subset $K$ of the positive orthant, choosing $a$ and $b$ adapted to $K$ gives a bound for the number of positive solutions in $K$ which is at most the bound of \cite{BDF}, together with a parity statement for the positive solutions in an explicit region containing $K$ (Corollary \ref{cor:compact}). This bound is not a rule of signs: it refines the bound of \cite{BDF} inside each set of coefficient matrices with a given pattern of signs of maximal minors, and it describes explicit semi-algebraic regions of coefficient matrices with fewer positive solutions.

\subsection*{Organization}
In Section \ref{known} we recall the generalized Budan--Fourier theorem of \cite{C}. In Section \ref{sec:wronskian} we prove the Budan--Fourier theorem with Wronskians, the renumbering result for sequences satisfying Descartes' rule of signs, and we study in detail the functions arising from circuits (Subsection \ref{ssec:setting3}) and negative powers of linear functions (Subsection \ref{ssec:powers}). In Section \ref{sec:circuits} we recall the results of \cite{BDF} on systems supported on circuits and prove Theorem \ref{thm:localized}.

\subsection*{Use of AI}
An AI assistant (Claude Opus 5.5, Anthropic) helped the author in the comparison with Descartes' rule of signs, namely in the analysis of the endpoints and in the proofs of Propositions \ref{prop:endpoints} and \ref{prop:BFleqD} and Corollary \ref{cor:VleqD}, and in the editing of the text. It also helped in obtaining Theorem \ref{thm:powers} on negative powers of linear functions. The author takes full responsibility for the contents of this paper.

\subsection*{Data availability}
Data sharing is not applicable to this article, as no datasets were generated or analysed during the current study.

\subsection*{Conflict of interest}
The author has no conflicts of interest to declare that are relevant to the content of this article.

\section{Generalizations of the Budan--Fourier theorem}
\label{known}

We recall notations and results from \cite{C} that we will need.
Let $I$ be a non-empty open interval of $\R$. Consider the ring $\Ci(I)$ of infinitely differentiable real-valued functions on $I$.
Let $\mathbf{f}=(f_i)_{i \in \N}$ be any sequence of positive functions in $\Ci(I)$. For $i \in \N$ define the operator $D_{i}$ on $\Ci(I)$
by $D_{i}(g)=\left(\frac{g}{f_i} \right)'$. Let $d$ be a positive integer. A function $g \in \Ci(I)$ has $\mathbf{f}$-degree $d$ if

$$(D_{0}D_{1} \cdots D_{d})(g)=0 \quad \text{and} \quad (D_{1} \cdots D_{d})(g) \neq 0.$$
A function $g \in \Ci(I)$ has $\mathbf{f}$-degree $0$ if $D_{0}(g)=0$ and $g \neq 0$.

If $g$ has $\mathbf{f}$-degree $d$, define $\operatorname{Der}_\mathbf{f}(g)=(g_{0},\ldots,g_{d})$ where $g_{d}=g$ and $g_{i-1}=D_{i}(g_{i})=\left(\frac{g_{i}}{f_i} \right)'$ for $i=d,d-1,\ldots,1$.  Thus, $g_{d}=g$, $g_{d-1}=D_d(g)=\left(\frac{g}{f_{d}} \right)'$ and $g_{d-2}=\left(\frac{g_{d-1}}{f_{d-1}} \right)'=D_{d-1}D_d(g), \ldots,
g_{0}=(D_{1}\cdots D_{d})(g)$.
The function $g_{i}$ has $\mathbf{f}$-degree $i$. The $\mathbf{f}$-multiplicity of a root $a$ of $g$ is defined as the index $m$ such that
$g_{d}(a)=\cdots=g_{d-m+1}(a)=0$ and $g_{d-m}(a) \neq 0$.

\begin{example}(\cite[Example 1]{C})
When all $f_i$ are equal to the constant function $1$, we recover the usual derivatives, degree and multiplicity.
\end{example}

Define the sign variation $\sgnvar(c_0, \dots, c_r)$ of a finite sequence $(c_0, \dots, c_r)$ of real numbers
as the number of pairs $(i,j)$ of integers between $0$ and $r$ which satisfy that $i<j$, $c_i \cdot c_j <0$ and $c_\ell=0$ for any index
$\ell$ with $i<\ell<j$. If $g \in \Ci(I)$ has $\mathbf{f}$-degree $d$ with $\operatorname{Der}_\mathbf{f}(g)=(g_{0},\ldots,g_{d})$ then for
any $x \in I$ define $\operatorname{Der}_\mathbf{f}(g; x)=\sgnvar(g_{0}(x),\ldots,g_{d}(x))$. For $x,y \in I$ such that $x<y$ define
$\operatorname{Der}_\mathbf{f}(g; (x,y))=\operatorname{Der}_\mathbf{f}(g; x)-\operatorname{Der}_\mathbf{f}(g; y)$.

\begin{example}(\cite[Example 2]{C})
\label{ex:monomials}
Consider a generalized univariate polynomial
$$g(x)=\sum_{i=0}^{d}c_{i}x^{a_{i}}$$
where $c_0,\ldots,c_{d},a_{0},\ldots,a_{d}$ are real numbers and $a_{0} < a_{1} < \cdots < a_{d}$.
Taking for $\mathbf{f}$ the sequence $(x^{a_{d}-a_{d-1}-1},
x^{a_{d-1}-a_{d-2}-1},\ldots, x^{a_{1}-a_{0}-1}, x^{a_0},1,1,\ldots)$ we get that $g$ has $\mathbf{f}$-degree $d$. Moreover,
$\operatorname{Der}_\mathbf{f}(g; x)=\sgnvar(c_{0},\ldots,c_{d})$ for any sufficiently small positive real number $x$ and $\operatorname{Der}_\mathbf{f}(g; x)=0$ for any sufficiently large positive real number $x$.
\end{example}

For $a,b \in I$ with $a<b$ denote by $n_\mathbf{f}(g;(a,b])$ the number of roots of $g$ in the half-open interval $(a,b]$ counted with $\mathbf{f}$-multiplicities.

The following result specializes to the classical Budan--Fourier theorem by taking for $\mathbf{f}$ the sequence where all $f_i$ are equal to the constant function $1$.

 \begin{thm} (\cite[Theorem 1]{C})
 \label{th:BF}
 Let $g \in \Ci(I)$ be a function having $\mathbf{f}$-degree $d$. For any $a,b \in I$ with $a<b$ we have
 \begin{enumerate}
 \item $n_\mathbf{f}(g;(a,b]) \leq \operatorname{Der}_\mathbf{f}(g; (a,b))$
  \item $\operatorname{Der}_\mathbf{f}(g; (a,b))-n_\mathbf{f}(g;(a,b])$ is even.
  \end{enumerate}
 \end{thm}

\section{A Budan--Fourier theorem with Wronskians}
\label{sec:wronskian}

We first observe that the $\mathbf{f}$-multiplicity defined in Section \ref{known} is in fact equal to the classical multiplicity.

\begin{lemma}
\label{f=classical}
The $\mathbf{f}$-multiplicity of a root $a$ of $g$ is equal to its classical multiplicity.
\end{lemma}
\begin{proof}
Set $u_i=\frac{1}{f_i}$. Then all $u_i$ are positive on $I$ and
$g_{i-1}=D_{i}(g_{i})=\left(u_ig_{i}\right)'$ for $i=d,d-1,\ldots,1$.
We have $g_{d-1}=(u_dg)'=u_d'g+u_dg'$ and $g_{d-2}=(u_{d-1}g_{d-1})'=(u_{d-1}'u_d'+u_{d-1}u_d^{''})g+(u_{d-1}'u_d+2u_{d-1}u_d')g'+u_{d-1}u_dg^{''}$.  It follows then by induction on $k \geq 1$
that
$$g_{d-k}=\sum_{\ell=0}^k T_{k\ell} \cdot g^{(\ell)},$$
where $T_{k \ell}$ is a sum of products of derivatives up to order $k$ of functions taken among $u_d,u_{d-1},\ldots,u_{d-k+1}$ and $T_{kk}=u_d \cdot u_{d-1} \cdots u_{d-k+1}$. Note that $T_{kk}$ does not vanish on $I$.
Thus, for any integer $m\geq 0$ and any $a \in I$ there exists an invertible triangular matrix $T(a)$ such that the column vector $(g_{d-k}(a))_{k=0,\ldots,m}$ is equal to the product of $T(a)$ with the column vector $(g^{(d-k)}(a))_{k=0,\ldots,m}$.
Since $T(a)$ is invertible we get $g(a)=g'(a)=\cdots=g^{(m)}(a)=0$ if and only if $g_{d}(a)=\cdots=g_{d-m}(a)=0$.
\end{proof}

\begin{ex}
We use Example \ref{ex:monomials}, taking for $g$ the polynomial $x(x-3)^{2}=9x-6x^{2}+x^{3}$ and for $\mathbf{f}$ the sequence $(f_0,f_1,f_2,f_3,\ldots)=(1,1,x,1,1,\ldots)$.
Then $g_2=g$, $g_1=(g/f_2)'=(9-6x+x^2)'=-6+2x$ and $g_0=(g_1/f_1)'=g_1'=2$ so $g$ has $\mathbf{f}$-degree $2$. A simple computation shows that
 both the classical multiplicity and the $\mathbf{f}$-multiplicity of $3$ as a root of $g$ are equal to $2$, as predicted by Lemma \ref{f=classical}. Note however that the classical degree and
 the $\mathbf{f}$-degree of $g$ are different.
\end{ex}

Consider analytic functions $h_1,\ldots,h_{s}: I \rightarrow \R$ where $I \subset \R$ is a non-empty open interval.
The Wronskian of $h_1,\ldots,h_s$ is the following determinant
$$
W(h_1,\ldots,h_s)= \det
\left(
\begin{array}{cccc}
h_1 & h_2 & \ldots & h_s \\
h_1'& h_2' & \ldots & h_s'\\
\vdots & \vdots & \ldots & \vdots \\
h_1^{(s-1)} & h_2^{(s-1)} & \ldots & h_s^{(s-1)}
\end{array}
\right).$$

Let $g:I \rightarrow \R$ be any analytic function. Define
\begin{equation}
\label{seq H}
H_{k}(g)=W(h_{1},\ldots,h_{k-1},g), \quad k=1,\ldots,s.
\end{equation}
Here $H_{1}(g)=g$ by convention.
Set $H(g)=(H_{1}(g),\ldots,H_{s}(g))$. Thus for $x \in I$ we have $H(g)(x)=(H_{1}(g)(x),\ldots,H_{s}(g)(x))$.
For $a<b$ in $I$, denote by $n(g;(a,b])$ the number of roots of $g$ in $(a,b]$ counted with multiplicities, and for an interval $J \subset I$ denote similarly by $n(g;J)$ the number of roots of $g$ in $J$.

\begin{thm}
\label{application}
Let $h_{1},\ldots,h_{s}$ be analytic functions on $I$ such that
$W(h_{1},\ldots,h_{k}) >0$ on $I$ for $k=1,\ldots,s$ (in particular $h_1=W(h_1)$ is positive on $I$). Let $a,b \in I$ with $a<b$. Let $c_{1},\ldots,c_{s}$ be real numbers with $c_{s} \neq 0$
and set $g=\sum_{i=1}^{s}c_ih_{i}$. Then
$$n(g;(a,b]) \leq \sgnvar H(g)(a)-\sgnvar H(g)(b)$$
and the difference $\sgnvar H(g)(a)-\sgnvar H(g)(b)-n(g;(a,b])$ is an even integer.
\end{thm}

 \begin{proof}
We introduce
 the sequence $\mathbf{f}=(f_{i})_{i \in \N}$ defined by
$ f_{i}=1$ for $i \geq s$, $f_{s-1}=h_{1}, f_{s-2}=\frac{W(h_1,h_2)}{h_{1}^{2}}$
 and
 \begin{equation}
 \label{seq f}
 f_{s-t}=\frac{W(h_{1},\ldots,h_{t-2}) \cdot W(h_{1},\ldots,h_{t-2}, h_{t-1}, h_{t})}{W(h_{1},\ldots,h_{t-1})^{2} }, \quad t=3,\ldots,s.
 \end{equation}
Set
$g_{s-1}=g$ and $g_{i-1}=D_{i}(g_{i})=\left(\frac{g_{i}}{f_i} \right)'$ for $i=s-1,s-2,\ldots,0$.
We compute that
$g_{s-2}=D_{s-1}(g)=\left(\frac{g}{f_{s-1}} \right)'=\left(\frac{g}{h_{1}} \right)'=\sum_{i=2}^s c_{i }\frac{W(h_1,h_{i})}{h_{1}^2}$. By \cite{P-S} (see page 108), for any analytic functions $\psi_{1},\ldots,\psi_{m}$
on $I$ we have
$$\left(
\frac
{W(\psi_{1},\ldots,\psi_{m-2},\psi_m)}
{W(\psi_{1},\ldots,\psi_{m-2},\psi_{m-1})}\right)^{'}=
\frac
{W(\psi_{1},\ldots,\psi_{m-2})W(\psi_{1},\ldots,\psi_{m})}{W(\psi_{1},\ldots,\psi_{m-1})^2}.
$$
Using this equality, we show by induction that
\begin{equation}
\label{used}
g_{s-t}=\frac{W(h_{1},\ldots,h_{t-2})}{W(h_{1},\ldots,h_{t-1})^{2} } \cdot \sum_{i=t}^s c_{i} W(h_{1},\ldots,h_{t-1},h_{i})\quad  t=3,\ldots,s.
\end{equation}
In particular, for $t=s$ we get $g_{0}=(D_{1} \cdots D_{s-1})(g)=c_{s} f_{0}$. But
$c_{s} \neq 0$, and all the $f_{i}$ are positive on $I$ since they are $1$, $h_1=W(h_1)$ or quotients of products of the Wronskians $W(h_1,\ldots,h_k)$, which are positive by assumption.
Thus $g_{0} \neq 0$. Finally, we have $(D_{0} \cdots D_{s-1})(g)=\left(\frac{c_{s} f_{0} }{f_{0}} \right)'=0$, and we conclude that $g$ has $\mathbf{f}$-degree $s-1$.
Applying Theorem \ref{th:BF} yields that $n_\mathbf{f}(g;(a,b]) \leq \operatorname{Der}_\mathbf{f}(g; (a,b))$
and that the difference between these numbers is an even integer. We have $\sum_{i=t}^s c_{i} W(h_{1},\ldots,h_{t-1},h_{i})=H_t(g)$ and
$\frac{W(h_{1},\ldots,h_{t-2})}{W(h_{1},\ldots,h_{t-1})^{2} } $ is positive on $I$. Thus, using \eqref{used},  we obtain $\operatorname{Der}_\mathbf{f}(g; (a,b))=\sgnvar H(g)(a)-\sgnvar H(g)(b)$.
Finally, we get $n_\mathbf{f}(g;(a,b]) =n(g;(a,b])$ by Lemma \ref{f=classical}.
 \end{proof}

\begin{rem}
\label{rem:Hk}
As already noticed in the proof of Theorem \ref{application}, if $g=\sum_{i=1}^s c_ih_i$ then
$H_k(g)= \sum_{i=k}^s c_{i} W(h_{1},\ldots,h_{k-1},h_{i})$, $k=1, \ldots,s$.
\end{rem}

We now recall the characterization of sequences of functions satisfying Descartes' rule of signs.

A sequence $(h_1,h_2,\ldots,h_s)$ of real-valued analytic functions defined
on an open interval $I \subset \R$ satisfies Descartes' rule of signs on $I$ if for any
sequence $c=(c_1,c_2,\ldots,c_s)$ of real numbers, the number of roots
of $c_1h_1+c_2h_2+\cdots+c_sh_s$ in $I$ counted with multiplicity never exceeds $\sgnvar(c)$.
The classical univariate Descartes' rule of signs asserts that monomial bases
satisfy Descartes' rule of signs on the open interval $(0,+\infty)$.
Observe that $(h_1,h_2,\ldots,h_s)$ satisfies Descartes' rule of signs
on $I$ if and only if $(h_s,h_{s-1},\ldots,h_1)$ has this property.

\begin{prop}[{\cite[Part V, items 87 and 90]{P-S}}]\label{P:Descartes}
A sequence of analytic functions $h_1,\ldots,h_s$ satisfies Descartes' rule of signs on $I \subset \R$
if and only if for any collection of integers $1 \leq j_1 <j_2 < \cdots < j_\ell \leq s$ we have
\begin{enumerate}
\item \label{E:cond1}
$\quad W(h_{j_1},\ldots,h_{j_\ell}) \neq 0$ on $I$,
\end{enumerate}
and for any collections of integers $1 \leq j_1 <j_2 < \cdots < j_\ell \leq s$
and $1 \leq j_1' <j_2' < \cdots < j_\ell' \leq s$ of the same size, we have
\begin{enumerate}
\setcounter{enumi}{1}
\item \label{E:cond2}
$\quad W(h_{j_1},\ldots,h_{j_\ell}) \cdot W(h_{j_1'},\ldots,h_{j_\ell'}) >0$ on $I$.\hfill\qed
\end{enumerate}
\end{prop}
If a sequence of analytic functions $h_1,\ldots,h_s$ satisfies Descartes' rule of signs on $I$, then $h_1, \dots, h_s$
do not vanish, have the same sign on $I$ and no two of them are proportional.

We now show that Theorem \ref{application} applies, after a suitable renumbering, to any sequence of functions satisfying Descartes' rule of signs (see \cite[Part V]{P-S} for instances of such sequences).

\begin{prop}
\label{prop:renumbering}
Let $h_1,\ldots,h_s$ be analytic functions satisfying Descartes' rule of signs on an open interval $I$.
Then there exist a permutation $\sigma$ of $\{1,\ldots,s\}$ and $e \in \{-1,1\}$ such that the functions
$\tilde h_i=e\, h_{\sigma(i)}$, $i=1,\ldots,s$, satisfy
$$W(\tilde h_1,\ldots,\tilde h_k)>0 \quad \text{on} \; I, \quad k=1,\ldots,s.$$
In particular $\tilde h_1,\ldots,\tilde h_s$ satisfy the hypotheses of Theorem \ref{application}.
\end{prop}

\begin{proof}
By Proposition \ref{P:Descartes}, for each $k=1,\ldots,s$ the Wronskians $W(h_{j_1},\ldots,h_{j_k})$ with $j_1<\cdots<j_k$ do not vanish on $I$, and they all have the same constant sign $\epsilon_k \in \{-1,1\}$ on $I$ (each of them is continuous and does not vanish on the interval $I$, and any two of them have the same sign by condition \eqref{E:cond2}).

Let $\sigma$ be any permutation of $\{1,\ldots,s\}$. For $k=1,\ldots,s$, let $\pi_k$ be the permutation which sorts $(\sigma(1),\ldots,\sigma(k))$ in increasing order. Since a Wronskian is a determinant whose columns are indexed by the functions, we get
$$W(h_{\sigma(1)},\ldots,h_{\sigma(k)})=\sgn(\pi_k)\,\epsilon_k \cdot |W(h_{\sigma(1)},\ldots,h_{\sigma(k)})| \quad \text{on} \; I.$$
The sign of $\pi_k$ is $(-1)^{\mathrm{inv}_k}$, where $\mathrm{inv}_k$ is the number of inversions of the sequence $(\sigma(1),\ldots,\sigma(k))$, that is, of pairs $i<j\leq k$ with $\sigma(i)>\sigma(j)$. Setting
$$m_j=\#\{i<j \, | \, \sigma(i)>\sigma(j)\}, \quad j=1,\ldots,s,$$
we have $\mathrm{inv}_k=m_1+\cdots+m_k$. Moreover, since $W(eh_{\sigma(1)},\ldots,eh_{\sigma(k)})=e^kW(h_{\sigma(1)},\ldots,h_{\sigma(k)})$, the sign of $W(\tilde h_1,\ldots,\tilde h_k)$ is
$$e^k\, (-1)^{m_1+\cdots+m_k}\, \epsilon_k.$$
Take $e=\epsilon_1$ and set $\tau_k=\epsilon_1^k\,\epsilon_k$ for $k=1,\ldots,s$, so that $\tau_1=1$. We have to find $\sigma$ such that $(-1)^{m_1+\cdots+m_k}=\tau_k$ for $k=1,\ldots,s$, that is, such that $m_1=0$ (which always holds) and
\begin{equation}
\label{eq:lehmer}
(-1)^{m_k}=\tau_{k-1}\tau_k, \quad k=2,\ldots,s.
\end{equation}
Define $m_1=0$ and, for $k=2,\ldots,s$, $m_k=0$ if $\tau_{k-1}\tau_k=1$ and $m_k=1$ otherwise. Then $0 \leq m_k \leq k-1$ for $k=1,\ldots,s$. It is well known that the map which associates to a permutation $\sigma$ of $\{1,\ldots,s\}$ its sequence $(m_1,\ldots,m_s)$ is a bijection onto the set of sequences of integers with $0 \leq m_k \leq k-1$. For the reader's convenience, here is the construction: $\sigma(s)$ is the unique value having exactly $m_s$ values larger than it, that is, $\sigma(s)=s-m_s$; then $\sigma(s-1)$ is the element of $\{1,\ldots,s\}\setminus\{\sigma(s)\}$ having exactly $m_{s-1}$ larger elements in this set, and so on. The permutation $\sigma$ obtained in this way satisfies \eqref{eq:lehmer}, which concludes the proof.
\end{proof}

\begin{rem}
\label{rem:renumbering}
\begin{enumerate}
\item In general, $\sigma$ cannot be taken to be the identity or the reversal $i \mapsto s+1-i$: for $s=3$, the product $\epsilon_1\epsilon_2\epsilon_3$ is invariant under reversing the order and under changing all signs (which multiply $\epsilon_k$ by $(-1)^{k(k-1)/2}$ and $(-1)^k$ respectively), and it can be equal to $-1$.
\item Since Descartes' rule of signs depends on the ordering of the functions, the sequence $\tilde h_1,\ldots,\tilde h_s$ does not satisfy Descartes' rule of signs in general. However, Theorem \ref{application} applies to $\tilde h_1,\ldots,\tilde h_s$. Note that, to apply Theorem \ref{application} to $g=\sum_{i} c_i h_i$, we need the coefficient of the last function $\tilde h_s=e\,h_{\sigma(s)}$ to be non-zero. If $c_{\sigma(s)}=0$, it suffices to apply Proposition \ref{prop:renumbering} to the subsequence of functions $h_j$ with $c_j \neq 0$, which still satisfies Descartes' rule of signs on $I$.
\item The converse of Proposition \ref{prop:renumbering} does not hold: the positivity of the Wronskians $W(\tilde h_1,\ldots,\tilde h_k)$ for $k=1,\ldots,s$ does not imply condition \eqref{E:cond2} of Proposition \ref{P:Descartes} for all subsequences. For instance, $h_1=1$, $h_2=t$ and $h_3=(t-\frac{3}{2})^2+1$ are positive on $I=(1,2)$ with $W(h_1)=W(h_1,h_2)=1$ and $W(h_1,h_2,h_3)=2$, but $W(h_1,h_3)=2t-3$ changes sign on $I$.
\end{enumerate}
\end{rem}

\begin{ex}
\label{ex:monomialsW}
Consider functions $h_i: (0,+\infty) \rightarrow \R$ defined by $h_i=x^{a_i}$ for $i=1,\ldots,s$, where the exponents $a_i$ are real numbers such that $a_1<a_2< \cdots<a_s$.
Then $W(h_1,\ldots,h_k)$ factors as $x^{a_1+\cdots+a_k-\frac{k(k-1)}{2}}$ times the determinant
of the matrix with $(i,j)$ entry the falling factorial $a_j(a_j-1)\cdots(a_j-i+2)$ for $i,j=1,\ldots,k$. This determinant is equal to the Vandermonde determinant of $a_1,\ldots,a_k$, and thus
\begin{equation}
\label{Wrons}
W(h_1,\ldots,h_k)=
x^{a_1+\cdots+a_k-\frac{k(k-1)}{2}}\cdot \prod_{1 \leq i < j \leq k} (a_j-a_i)
.
\end{equation}
In particular $W(h_1,\ldots,h_k)$ is positive on $(0,+\infty)$ for $k=1,\ldots,s$ and we may apply Theorem \ref{application} on any interval $(a,b) \subset (0,+\infty)$.

Let $g=\sum_{i=1}^s c_ih_i$ with $c_s \neq 0$, and fix $t \in \{1,\ldots,s\}$. By Remark \ref{rem:Hk},
$$H_t(g)=\sum_{i=t}^s c_{i} W(h_{1},\ldots,h_{t-1},h_{i}),$$
and by \eqref{Wrons}, for $t \leq i<j \leq s$ the quotient
$$\frac{W(h_1,\ldots,h_{t-1},h_j)}{W(h_1,\ldots,h_{t-1},h_i)}$$
is equal to $x^{a_j-a_i}$ up to a positive multiplicative constant. Since $a_j-a_i >0$, for $x>0$ large enough the sign of
$H_t(g)(x)$ equals that of $c_s$, and thus $\sgnvar H(g)(x)=0$ for $x>0$ large enough.
Similarly, for $x>0$ small enough the sign of $H_t(g)(x)$ equals the sign of the first non-zero coefficient among $c_t,c_{t+1},\ldots,c_s$. Since replacing each zero term of $(c_1,\ldots,c_s)$ by the next non-zero term does not change the number of sign changes, we get $\sgnvar H(g)(x)=\sgnvar(c_1,\ldots,c_s)$ for $x>0$ small enough.
Therefore, we recover Descartes' rule of signs when applying Theorem \ref{application} to $(0,+\infty)$.

\end{ex}

\begin{ex}
\label{firstD}
Let $p_i=a_i+b_iy \in \R[y]$, $i=1,\ldots,s$. Assume that the determinants $a_ib_j-a_jb_i$ are negative whenever $i<j$. We also assume that the interval
$I=\{y \in \R \, | \, p_1(y) >0,\ldots,p_s(y)>0\}$ is non-empty. Consider the functions $h_i=\frac{1}{p_i}$ for $i=1,\ldots,s$. These functions play a central role in \cite{BD}, which contains the following Wronskian computation.
We have
$$
(1/p_j)^{(i-1)} =(-1)^{i-1} (i-1)! \cdot\frac{b_j^{i-1}}{p_j^i}.
$$
Thus, a Vandermonde determinant computation gives
\begin{align*}
W(h_1,\ldots,h_{k}) & = \kappa_k\cdot \frac{1}{p_1\cdots p_{k}} \cdot \det
\Big( \Big(\frac{b_j}{p_j}\Big)^{i-1}\Big)_{1 \leq i,j \leq k} \\
& = \kappa_k\cdot \frac{1}{p_1\cdots p_k} \cdot \prod_{1 \leq r< t \leq k} \Big(\frac{b_t}{p_{t}}-\frac{b_r}{p_{r}}\Big),
\end{align*}
where $\kappa_k\, = \, \prod_{i=1}^{k}  (-1)^{i-1} (i-1)! $.
Note that $\frac{b_t}{p_t}-\frac{b_r}{p_{r}}=\frac{a_rb_t-b_ra_t}{p_tp_r}$ is negative on $I$ if $r<t$.
It follows that $W(h_1,\ldots,h_k)$ does not vanish on $I$ and its sign on $I$ is equal to the sign of $\prod_{i=1}^{k}  (-1)^{i-1} \cdot (-1)^{k(k-1)/2}=1$.
Thus all Wronskians $W(h_1,\ldots,h_k)$ are positive on $I$ and we can apply Theorem \ref{application}.
\end{ex}

\begin{ex}
\label{secD}
Let $p_1,\ldots,p_s$ and $I$ be as in Example \ref{firstD}.
Set $h_{i}=\frac{b_i}{p_i}-\frac{b_{i+1}}{p_{i+1}}$ for $i=1,\ldots, s-1$. These functions have been used in \cite{BDF}, where the Wronskian
$W(h_{1},\ldots,h_k)$ has been computed:
we have
$$
\bigg(\frac{b_j}{p_j}\bigg)^{(i-1)} =(-1)^{i-1}\, (i-1)!  \bigg(\frac{b_j}{p_j}\bigg)^{i},
$$
and using a Vandermonde determinant we get

\begin{align*}
W(h_{1},\ldots,h_{k})& =
(-1)^k \kappa_k \cdot \det
\left( \Big(\frac{b_j}{p_j}\Big)^{i-1}
\right)_{1\leq i,j \leq k+1} \\
& = \displaystyle(-1)^{k} \kappa_k \cdot \prod_{1 \leq r < t \leq k+1} \Big(\frac{b_t}{p_t}-\frac{b_r}{p_r}\Big),
\end{align*}
where $\kappa_k$ is defined in Example \ref{firstD}.
Again, we have $\frac{b_t}{p_{t}}-\frac{b_r}{p_{r}}=\frac{a_rb_t-b_ra_t}{p_tp_r}$, so this function is negative on $I$ if $r<t$.
It follows that $W(h_{1},\ldots,h_{k})$ does not vanish on $I$ and its sign on $I$ is given by that of
$(-1)^{k}\kappa_k \cdot (-1)^{\frac{k(k+1)}{2}}$, which is positive.
So again we conclude that all $W(h_{1},\ldots,h_{k})$ are positive on $I$.
\end{ex}

It is shown in \cite{BDF}, using Proposition \ref{P:Descartes}, that the functions $h_1,\ldots,h_{s-1}$ of Example \ref{secD} satisfy Descartes' rule of signs on $I$.

\subsection{Linear combinations of the functions of Example \ref{secD}}
\label{ssec:setting3}
We now study in detail the bound given by Theorem \ref{application} for the functions of Example \ref{secD}. Throughout this subsection we fix linear functions
$p_j(y)=a_j+b_jy=\langle P_j,(1,y)\rangle$, where $P_j=(a_j,b_j) \in \R^2$, $j=1,\ldots,s$, with $s \geq 2$, and we assume
\begin{equation}
\label{hyp:L}
\tag{$\mathcal{L}$}
\begin{gathered}
\det(P_i,P_j)=a_ib_j-a_jb_i<0 \;\; \text{for } 1 \leq i<j \leq s, \\
I=\{y \in \R \, | \, p_1(y)>0,\ldots,p_s(y)>0\} \neq \emptyset.
\end{gathered}
\end{equation}
We use the following notations:
\begin{itemize}
\item $N=s-1$, $q_j=\frac{b_j}{p_j}$ for $j=1,\ldots,s$, and $h_i=q_i-q_{i+1}$ for $i=1,\ldots,N$;
\item for $c=(c_1,\ldots,c_N) \in \R^N$, $g=\sum_{i=1}^N c_ih_i$, $g_2=\sum_{i=2}^N c_ih_i$, and $H(g)=(H_1(g),\ldots,H_N(g))$ is the sequence \eqref{seq H} associated with $h_1,\ldots,h_N$.
\end{itemize}
By Example \ref{secD}, $W(h_1,\ldots,h_k)>0$ on $I$ for $k=1,\ldots,N$. Since
$q_i-q_j=\frac{a_jb_i-a_ib_j}{p_ip_j}$, we have
$q_1>q_2>\cdots>q_s$ on $I$; in particular all $h_i$ are positive on $I$. Moreover
\begin{equation}
\label{derq}
q_j^{(m)}=(-1)^m\, m!\, q_j^{m+1}, \quad m \geq 0.
\end{equation}

\smallskip

\noindent\emph{The bounded case.} The interval $I$ is bounded if and only if $b_1>0$ and $b_s<0$. Indeed, $I$ is an intersection of half-lines, bounded from below by the conditions $p_j>0$ with $b_j>0$ and from above by those with $b_j<0$ (the conditions with $b_j=0$ hold on $\R$ since $I \neq \emptyset$); thus $I$ is bounded if and only if some $b_j$ is positive and some $b_j$ is negative. Moreover, for $y_0 \in I$ the numbers $q_j(y_0)$ are strictly decreasing in $j$ and have the signs of the $b_j$, so the signs of $b_1,\ldots,b_s$ are non-increasing. Hence some $b_j$ is positive if and only if $b_1>0$, and some $b_j$ is negative if and only if $b_s<0$.

Assume that $I=(\alpha,\beta)$ is bounded. At an endpoint of $I$ some $p_j$ vanishes. At $\alpha$ the function $p_j$ which vanishes is increasing, so that $q_j \to +\infty$ when $y \to \alpha^+$; since $q_1$ is the largest of the $q_j$, this forces $p_1(\alpha)=0$.
Similarly $q_j\to -\infty$ at $\beta$ and thus $p_s(\beta)=0$. Since $a_ib_j-a_jb_i \neq 0$ for $i \neq j$, two distinct $p_j$ have distinct roots. Hence
$p_2(\alpha),\ldots,p_s(\alpha)>0$ and $p_1(\beta),\ldots,p_{s-1}(\beta)>0$.
In particular $\alpha$ belongs to $I_2=\{y \in \R \, | \, p_2(y)>0,\ldots,p_s(y)>0\}$, on which $h_2,\ldots,h_N$ satisfy the hypotheses of Example \ref{secD}. For $c \in \R^N$ we set
$H^{(2)}_k(g_2)=W(h_2,\ldots,h_k,g_2)$ for $k=1,\ldots,N-1$ (so $H^{(2)}_1(g_2)=g_2$ and $H^{(2)}_{N-1}(g_2)=c_N W(h_2,\ldots,h_N)$),
\begin{equation}
\label{eq:Valpha}
T_k(c)=H^{(2)}_k(g_2)(\alpha), \quad k=1,\ldots,N-1, \qquad \text{and} \qquad V_\alpha(c)=\sgnvar\big(c_1,T_1(c),\ldots,T_{N-1}(c)\big).
\end{equation}
Note that $T_{N-1}(c)$ has the sign of $c_N$.

\begin{prop}
\label{prop:endpoints}
Assume \eqref{hyp:L} and that $I=(\alpha,\beta)$ is bounded. Let $c \in \R^N$ with $c_1 \neq 0$, $c_N \neq 0$ and $T_k(c) \neq 0$ for $k=1,\ldots,N-1$ (see \eqref{eq:Valpha}). Then
\begin{enumerate}
\item $\sgnvar H(g)(y)=0$ for $y<\beta$ close enough to $\beta$,
\item $\sgnvar H(g)(y)=V_\alpha(c)$ for $y>\alpha$ close enough to $\alpha$.
\end{enumerate}
\end{prop}

\begin{proof}
By \eqref{derq}, if $\varphi=\sum_{j=1}^s v_j q_j$ then $\varphi^{(m)}=(-1)^m m! \sum_j v_j q_j^{m+1}$. Thus, for functions $\varphi_1,\ldots,\varphi_k$ in the linear span of $q_1,\ldots,q_s$, the Wronskian $W(\varphi_1,\ldots,\varphi_k)$ is multilinear in the columns $(q_j,-q_j^2,\ldots,(-1)^{k-1}(k-1)!\,q_j^k)^{T}$, and any term in which the same $q_j$ appears in two columns vanishes.

Let us first consider the endpoint $\beta$. When $y \to \beta^-$, the function $q_s$ tends to $-\infty$ while $q_1,\ldots,q_{s-1}$ remain bounded. For $k \leq N$, the function $q_s$ appears in $H_k(g)=W(h_1,\ldots,h_{k-1},g)$ only through the term $c_N h_N=c_N(q_N-q_s)$ of $g$ (if $k=N$ then $h_{N-1}=q_{N-1}-q_N$ does not contain $q_s$ either). Expanding with respect to the last column, the dominant term is obtained by taking the entry $-c_N(-1)^{k-1}(k-1)!\,q_s^k$ of the last row, and we get
$$
H_k(g)(y)= c_N\, (k-1)!\, |q_s(y)|^k \, W(h_1,\ldots,h_{k-1})(y) \,(1+o(1)), \quad y \to \beta^-.
$$
Since $\beta$ belongs to the interval $\{p_1>0,\ldots,p_{s-1}>0\}$, on which Example \ref{secD} applies to $h_1,\ldots,h_{N-1}$, we have $W(h_1,\ldots,h_{k-1})(\beta)>0$. Hence all $H_k(g)$ have the sign of $c_N$ near $\beta$, which proves the first item.

Consider now the endpoint $\alpha$. When $y \to \alpha^+$, the function $q_1$ tends to $+\infty$ while $q_2,\ldots,q_s$ remain bounded. First, $H_1(g)=g=c_1q_1+O(1)$ has the sign of $c_1$ near $\alpha$. For $k \geq 2$ we have
$H_k(g)=W(h_1,\ldots,h_{k-1},g-c_1h_1)=W(h_1,\ldots,h_{k-1},g_2)$, and $q_1$ appears only in the first column through $h_1=q_1-q_2$. Expanding with respect to the first column, the dominant term is given by the entry $(-1)^{k-1}(k-1)!\,q_1^k$ of the last row, whose cofactor is $(-1)^{k+1}$ times the Wronskian $W(h_2,\ldots,h_{k-1},g_2)$ (the latter is formed by the first $k-1$ rows). Thus
$$
H_k(g)(y)=(k-1)!\, q_1(y)^k\, H^{(2)}_{k-1}(g_2)(\alpha)\,(1+o(1)), \quad y \to \alpha^+,
$$
which proves the second item.
\end{proof}

We now show that the Budan--Fourier sign sequence at any point of $I$ never exceeds Descartes' bound. We do not assume anymore that $I$ is bounded. We first need an invariance property under projective changes of variable.

\begin{lemma}
\label{lem:projective}
Assume \eqref{hyp:L} and $0 \in I$. Let $\lambda \in \R$, and let $\tilde p_j$, $\tilde I$, $\tilde h_i$ be defined as $p_j$, $I$, $h_i$, starting from the vectors
$\tilde P_j=(a_j,b_j-\lambda a_j)$ instead of $P_j$. Then the $\tilde p_j$ satisfy \eqref{hyp:L}, $0 \in \tilde I$, and for any $c \in \R^N$, setting $\tilde g=\sum_{i=1}^Nc_i\tilde h_i$, we have
$$H(\tilde g)(0)=H(g)(0),$$
where $H(\tilde g)$ is the sequence \eqref{seq H} associated with $\tilde h_1,\ldots,\tilde h_N$.
\end{lemma}

\begin{proof}
Set $\varphi(t)=\frac{t}{1-\lambda t}$. We have $\tilde P_j=P_j\cdot M$ with $M=\begin{pmatrix}1 & -\lambda \\ 0 & 1\end{pmatrix}$, whence $\det(\tilde P_i,\tilde P_j)=\det(P_i,P_j)$ since $\det M=1$. Moreover $\tilde p_j(0)=a_j=p_j(0)>0$, so $0 \in \tilde I$.
On the open interval $U$ containing $0$ where $1-\lambda t>0$ and $\varphi(t) \in I$, we have $\tilde p_j(t)=(1-\lambda t)\, p_j(\varphi(t))$. Since $q_j=\frac{p_j'}{p_j}$ and $\tilde q_j=\frac{\tilde p_j'}{\tilde p_j}$, taking logarithmic derivatives we get
$\tilde q_j(t)=\frac{-\lambda}{1-\lambda t}+\varphi'(t)\, q_j(\varphi(t))$, and thus
$$\tilde h_i(t)=\varphi'(t)\, h_i(\varphi(t)), \qquad \tilde g(t)=\varphi'(t)\, g(\varphi(t)), \qquad t\in U.$$
For any smooth functions $f_1,\ldots,f_k$ and any smooth function $w$ we have the classical identities
$W(wf_1,\ldots,wf_k)=w^k\,W(f_1,\ldots,f_k)$ and
$W(f_1\circ\varphi,\ldots,f_k\circ\varphi)=(\varphi')^{\frac{k(k-1)}{2}}\, W(f_1,\ldots,f_k)\circ\varphi$
(the second one follows from the fact that $(f\circ \varphi)^{(m)}=(\varphi')^m f^{(m)}\circ \varphi+\sum_{r<m}u_{m,r}\, f^{(r)}\circ \varphi$, where the functions $u_{m,r}$ do not depend on $f$). It follows that
$$H_k(\tilde g)(t)=\varphi'(t)^{k+\frac{k(k-1)}{2}}\, H_k(g)(\varphi(t)), \quad t \in U, \; k=1,\ldots,N.$$
Since $\varphi(0)=0$ and $\varphi'(0)=1$, we get $H(\tilde g)(0)=H(g)(0)$.
\end{proof}

\begin{prop}
\label{prop:BFleqD}
Assume \eqref{hyp:L} ($I$ may be unbounded) and let $c \in \R^N$ with $c_N \neq 0$. Then for any $y \in I$ we have
$$\sgnvar H(g)(y) \leq \sgnvar(c_1,\ldots,c_N).$$
\end{prop}

\begin{proof}
We proceed by induction on $N$. If $N=1$ then $H(g)=(g)$ and $\sgnvar H(g)(y)=0$. Assume $N \geq 2$ and that the result holds for families of $N-1$ functions. Fix $y \in I$.

\emph{Step 1: reduction to the case where $I$ is bounded.} Replacing $p_j(t)$ by $p_j(t+y)$ (a translation, which changes neither the hypotheses nor the Wronskians), we may assume $y=0$. Then $a_j=p_j(0)>0$ for all $j$, and the numbers $\frac{b_j}{a_j}$ are pairwise distinct since $\det(P_i,P_j) \neq 0$ for $i\neq j$. Choose $\lambda$ strictly between the smallest and the largest of the numbers $\frac{b_j}{a_j}$. With the notations of Lemma \ref{lem:projective}, the slopes $b_j-\lambda a_j$ of the linear functions $\tilde p_j$ take both signs, so that $\tilde I$ is bounded. By Lemma \ref{lem:projective}, the new family satisfies the hypotheses of Example \ref{secD} on $\tilde I \ni 0$ and $H(\tilde g)(0)=H(g)(0)$. Therefore we may assume that $I=(\alpha,\beta)$ is bounded.

\emph{Step 2: reduction to generic coefficients.} For $i=1,\ldots,N$ let $\epsilon_i$ be the sign of $c_i$ if $c_i \neq 0$, and the sign of the first non-zero $c_j$ with $j>i$ otherwise (it exists since $c_N \neq 0$). Then $\sgnvar(\epsilon_1,\ldots,\epsilon_N)=\sgnvar(c_1,\ldots,c_N)$, and $c$ lies in the closure of the open set $\mathcal{U}=\{c' \in \R^N \, | \, \sgn(c'_i)=\epsilon_i, \; i=1,\ldots,N\}$. For $c' \in \mathcal{U}$ we have $\sgnvar(c')=\sgnvar(c)$.
Each entry of $H(g)(y)$ depends linearly on $c$. Hence, for $c' \in \mathcal{U}$ close enough to $c$, the non-zero entries of $H(g)(y)$ have the same signs as the corresponding entries of $H(g')(y)$, where $g'=\sum_i c'_ih_i$. Since removing terms from a sequence cannot increase its number of sign changes, we get $\sgnvar H(g)(y) \leq \sgnvar H(g')(y)$.
For $k=1,\ldots,N-1$, the number
$T_k(c')=\sum_{i=k+1}^N c'_i\, W(h_2,\ldots,h_k,h_i)(\alpha)$ (see \eqref{eq:Valpha}) is a linear form in $c'$ whose coefficient of $c'_{k+1}$ is $W(h_2,\ldots,h_{k+1})(\alpha)>0$. Thus, we can choose $c' \in \mathcal{U}$ arbitrarily close to $c$ such that $T_k(c') \neq 0$ for $k=1,\ldots,N-1$. It is thus enough to prove the result for $c'$, that is, we may assume that $c_1,\ldots,c_N$ are non-zero and that $T_k(c) \neq 0$ for $k=1,\ldots,N-1$.

\emph{Step 3: comparison with the left endpoint.} By Theorem \ref{application}, for $\alpha<y'<y$ we have $\sgnvar H(g)(y')-\sgnvar H(g)(y) \geq n(g;(y',y]) \geq 0$. Letting $y'$ tend to $\alpha$ and using Proposition \ref{prop:endpoints}, we get
$$\sgnvar H(g)(y) \leq V_\alpha(c)=\sgnvar(c_1,T_1,\ldots,T_{N-1}), \quad \text{where} \; T_k=T_k(c).$$

\emph{Step 4: induction.} The functions $h_2,\ldots,h_N$ satisfy the hypotheses of Example \ref{secD} on $I_2=\{p_2>0,\ldots,p_s>0\}$, and $\alpha \in I_2$. By the induction hypothesis applied to $g_2=\sum_{i=2}^Nc_ih_i$ at the point $\alpha$,
$$\sgnvar(T_1,\ldots,T_{N-1}) \leq \sgnvar(c_2,\ldots,c_N).$$
All the numbers $T_k$ and $c_i$ are non-zero, and $T_{N-1}=c_NW(h_2,\ldots,h_N)(\alpha)$ has the sign of $c_N$. Recall that the number of sign changes of a sequence of non-zero numbers has the parity of the number of sign changes between its first and last terms.
If $T_1$ and $c_2$ have the same sign, then
$\sgnvar(c_1,T_1,\ldots,T_{N-1})=\sgnvar(c_1,T_1)+\sgnvar(T_1,\ldots,T_{N-1}) \leq \sgnvar(c_1,c_2)+\sgnvar(c_2,\ldots,c_N)=\sgnvar(c)$.
If $T_1$ and $c_2$ have opposite signs, then $\sgnvar(T_1,\ldots,T_{N-1})$ and $\sgnvar(c_2,\ldots,c_N)$ have different parities, so that $\sgnvar(T_1,\ldots,T_{N-1}) \leq \sgnvar(c_2,\ldots,c_N)-1$, and thus
$\sgnvar(c_1,T_1,\ldots,T_{N-1}) \leq 1+\sgnvar(T_1,\ldots,T_{N-1}) \leq \sgnvar(c_2,\ldots,c_N) \leq \sgnvar(c)$.
In both cases $\sgnvar H(g)(y) \leq V_\alpha(c) \leq \sgnvar(c)$.
\end{proof}

\begin{cor}
\label{cor:VleqD}
Assume \eqref{hyp:L} and that $I=(\alpha,\beta)$ is bounded. Let $c \in \R^N$ with $c_1 \neq 0$, $c_N \neq 0$ and $T_k(c) \neq 0$ for $k=1,\ldots,N-1$. Then
\begin{equation}
\label{VleqD}
n(g;I) \leq V_\alpha(c) \leq \sgnvar(c_1,\ldots,c_N), \quad \text{and} \quad n(g;I) \equiv V_\alpha(c) \equiv \sgnvar(c_1,\ldots,c_N) \mod 2.
\end{equation}
Thus the bound $n(g;I) \leq V_\alpha(c)$ refines Descartes' rule of signs $n(g;I) \leq \sgnvar(c_1,\ldots,c_N)$ for the functions $h_1,\ldots,h_N$ (Proposition \ref{P:Descartes}).
\end{cor}

\begin{proof}
By Proposition \ref{prop:endpoints}, $\sgnvar H(g)(y)=V_\alpha(c)$ for $y>\alpha$ close to $\alpha$ and $\sgnvar H(g)(y)=0$ for $y<\beta$ close to $\beta$. Since $|g| \to +\infty$ at $\alpha$ and at $\beta$, the roots of $g$ lie in a compact subinterval of $I$, and Theorem \ref{application} applied to $(a,b]$ with $a$ close to $\alpha$ and $b$ close to $\beta$ gives $n(g;I) \leq V_\alpha(c)$ with an even difference. The inequality $V_\alpha(c) \leq \sgnvar(c_1,\ldots,c_N)$ follows from Proposition \ref{prop:BFleqD}. Finally, the first term of the sequence defining $V_\alpha(c)$ has the sign of $c_1$ and its last term has the sign of $c_N$; since the parity of the number of sign changes of a sequence of non-zero numbers only depends on the signs of its first and last terms, we get $V_\alpha(c) \equiv \sgnvar(c_1,\ldots,c_N) \mod 2$.
\end{proof}

\begin{rem}
\label{rem:comparisonDescartes}
Inequality \eqref{VleqD} can be strict and $V_\alpha(c)$ is not determined by the signs of $c_1,\ldots,c_N$. Contrary to Example \ref{ex:monomialsW} (where the functions $x^{a_i}$ have pairwise distinct orders at $0$), only $h_1$ is singular at $\alpha$, and the values at $\alpha$ of the other functions $h_i$ enter the bound.
For $N=2$ we get $V_\alpha(c)=\sgnvar(c_1,c_2)$. For $N=3$ we get
$$V_\alpha(c)=\sgnvar\big(c_1,\, c_2h_2(\alpha)+c_3h_3(\alpha),\, c_3\big).$$
For instance, if $c_1,c_3>0>c_2$ then Descartes' rule of signs gives $n(g;I) \leq 2$, while Corollary \ref{cor:VleqD} gives $n(g;I)=0$ when $c_3h_3(\alpha)>|c_2|h_2(\alpha)$.
\end{rem}
\begin{rem}
\label{rem:unbounded}
The boundedness assumption in Proposition \ref{prop:endpoints} and Corollary \ref{cor:VleqD} is not a real restriction. For systems supported on circuits, the Gale dual configuration can always be chosen so that $\Delta_P$ is bounded (Lemma \ref{lem:bounded}). For general linear functions satisfying \eqref{hyp:L} with $I$ unbounded, Step 1 of the proof of Proposition \ref{prop:BFleqD} (a translation such that $0 \in I$, so that all $a_j>0$, followed by Lemma \ref{lem:projective} with $\lambda$ strictly between the smallest and the largest of the numbers $b_j/a_j$) produces a bounded interval $\tilde I$. For $y \in I$, the number $1+\lambda y$ is a convex combination of $p_j(y)/a_j>0$ for two indices $j$, hence positive, so $\varphi(t)=\frac{t}{1-\lambda t}$ maps an open subinterval $U \subset \tilde I$ bijectively and increasingly onto $I$, and $\tilde g=\varphi'\cdot g\circ\varphi$ on $U$. Therefore $n(g;I)=n(\tilde g;U) \leq n(\tilde g;\tilde I)$, and Corollary \ref{cor:VleqD} applied to $\tilde g$ on $\tilde I$ gives an upper bound for $n(g;I)$.
\end{rem}

\subsection{Negative powers of linear functions}
\label{ssec:powers}
The method of Subsection \ref{ssec:setting3} applies to other families of functions satisfying Descartes' rule of signs. We consider here powers of the linear functions $p_j$. Throughout this subsection we assume \eqref{hyp:L}, we fix a real number $\gamma$, and we set
$$f_j=p_j^{\gamma}, \quad j=1,\ldots,s, \qquad D_{ij}=\det(P_i,P_j)=a_ib_j-a_jb_i,$$
so that $D_{ij}<0$ for $i<j$. For $m \geq 0$ we denote by $(\gamma)_m=\gamma(\gamma-1)\cdots(\gamma-m+1)$ the falling factorial, with $(\gamma)_0=1$, and we set $\Gamma_k=\prod_{m=0}^{k-1}(\gamma)_m$. For $c=(c_1,\ldots,c_s) \in \R^s$ we write $F=\sum_{j=1}^s c_jf_j$, and $H(F)=(H_1(F),\ldots,H_s(F))$ is the sequence \eqref{seq H} associated with $f_1,\ldots,f_s$.

\begin{lemma}
\label{lem:powersW}
For $1 \leq j_1<\cdots<j_k \leq s$ we have on $I$
$$W(f_{j_1},\ldots,f_{j_k})=\Gamma_k \cdot \prod_{r=1}^k p_{j_r}^{\gamma-k+1} \cdot \prod_{1 \leq r<t \leq k} D_{j_rj_t}.$$
In particular, if $\gamma \notin \{0,1,\ldots,s-2\}$, the sign of $W(f_{j_1},\ldots,f_{j_k})$ on $I$ only depends on $k$, so that $f_1,\ldots,f_s$ satisfy Descartes' rule of signs on $I$ (Proposition \ref{P:Descartes}). If $\gamma<0$, all these Wronskians are positive on $I$.
\end{lemma}

\begin{proof}
We have $f_j^{(m)}=(\gamma)_m\, b_j^m\, p_j^{\gamma-m}$. Factoring out $(\gamma)_m$ from the row of index $m \in \{0,\ldots,k-1\}$ and $p_{j}^{\gamma-k+1}$ from the column of $f_j$, we are left with the determinant of the matrix with entries $b_{j_r}^m\,p_{j_r}^{k-1-m}$, which is the homogeneous Vandermonde determinant
$\prod_{r<t}\big(p_{j_r}b_{j_t}-p_{j_t}b_{j_r}\big)=\prod_{r<t}D_{j_rj_t}$, since $p_ib_j-p_jb_i=a_ib_j-a_jb_i$. The sign of $\Gamma_k\prod_{r<t}D_{j_rj_t}$ is that of $\Gamma_k\cdot(-1)^{k(k-1)/2}$, which only depends on $k$, and $\Gamma_k \neq 0$ when $\gamma \notin \{0,\ldots,k-2\}$. If $\gamma<0$, then $(\gamma)_m$ has the sign of $(-1)^m$, so that $\Gamma_k$ has the sign of $(-1)^{k(k-1)/2}$ and the Wronskian is positive.
\end{proof}

For $\gamma<0$, the sequence $H(F)$ has a simple closed form. For $k=1,\ldots,s$ and $y \in I$ we set
\begin{equation}
\label{eq:sigma}
\sigma_k(y)=\sum_{i=k}^{s} c_i\, \delta_{k,i}\, p_i(y)^{\gamma-k+1}, \qquad \text{where} \quad \delta_{k,i}=\prod_{j=1}^{k-1}|D_{ji}|>0,
\end{equation}
and $\sigma(y)=(\sigma_1(y),\ldots,\sigma_s(y))$. Note that $\sigma_1=F$.

\begin{prop}
\label{prop:sigma}
Assume $\gamma<0$. For $k=1,\ldots,s$ we have on $I$
$$H_k(F)=|\Gamma_k| \cdot \prod_{1 \leq j<j' \leq k-1}|D_{jj'}| \cdot \prod_{j=1}^{k-1}p_j^{\gamma-k+1} \cdot \sigma_k.$$
Consequently, if $c_s \neq 0$, then for any $a<b$ in $I$
$$n(F;(a,b]) \leq \sgnvar\sigma(a)-\sgnvar\sigma(b),$$
and the difference is an even integer.
\end{prop}

\begin{proof}
By Remark \ref{rem:Hk}, $H_k(F)=\sum_{i=k}^s c_i W(f_1,\ldots,f_{k-1},f_i)$, and by Lemma \ref{lem:powersW},
$$W(f_1,\ldots,f_{k-1},f_i)=\Gamma_k \prod_{j=1}^{k-1}p_j^{\gamma-k+1} \cdot p_i^{\gamma-k+1} \prod_{1 \leq j<j' \leq k-1}D_{jj'} \cdot \prod_{j=1}^{k-1}D_{ji}.$$
Since $D_{jj'}<0$ for $j<j'$ and $\Gamma_k$ has the sign of $(-1)^{k(k-1)/2}$, the sign of $\Gamma_k\prod_{j<j'<k}D_{jj'}\prod_{j<k}D_{ji}$ is $(-1)^{k(k-1)/2+(k-1)(k-2)/2+(k-1)}=(-1)^{k(k-1)}=1$. This gives the formula for $H_k(F)$. The second assertion follows from Theorem \ref{application} and Lemma \ref{lem:powersW}, since the factor in front of $\sigma_k$ is positive.
\end{proof}

\begin{thm}
\label{thm:powers}
Assume \eqref{hyp:L} and $\gamma<0$, and let $c \in \R^s$ with $c_s \neq 0$.
\begin{enumerate}
\item For any $y \in I$ we have $\sgnvar\sigma(y) \leq \sgnvar(c_1,\ldots,c_s)$.
\item Assume moreover that $I=(\alpha,\beta)$ is bounded, that $c_1 \neq 0$, and that $\sigma_k(\alpha) \neq 0$ for $k=2,\ldots,s$ (note that $\sigma_k$, $k \geq 2$, does not involve $p_1$ and is defined at $\alpha$). Set
$$V^\gamma_\alpha(c)=\sgnvar\big(c_1,\sigma_2(\alpha),\ldots,\sigma_s(\alpha)\big).$$
Then $\sgnvar\sigma(y)=0$ for $y$ close to $\beta$, $\sgnvar\sigma(y)=V^\gamma_\alpha(c)$ for $y$ close to $\alpha$, and
$$n(F;I) \leq V^\gamma_\alpha(c) \leq \sgnvar(c_1,\ldots,c_s), \qquad n(F;I) \equiv V^\gamma_\alpha(c) \equiv \sgnvar(c_1,\ldots,c_s) \mod 2.$$
Moreover $\sigma_2(\alpha)=b_1\sum_{i=2}^s c_i\,p_i(\alpha)^{\gamma}$.
\end{enumerate}
\end{thm}

\begin{proof}
We first prove the second item, assuming the first one. Recall from Subsection \ref{ssec:setting3} that $p_1(\alpha)=0$, $p_s(\beta)=0$, $p_2(\alpha),\ldots,p_s(\alpha)>0$ and $p_1(\beta),\ldots,p_{s-1}(\beta)>0$. When $y \to \beta^-$, each $\sigma_k$ contains the term $c_s\delta_{k,s}\,p_s^{\gamma-k+1}$, which tends to $\pm\infty$ since $\gamma-k+1<0$, while the other terms remain bounded; thus all $\sigma_k$ have the sign of $c_s$ near $\beta$. When $y \to \alpha^+$, $\sigma_1=F$ has the sign of $c_1$ since $p_1^\gamma \to +\infty$, while $\sigma_k$, $k \geq 2$, tends to $\sigma_k(\alpha)$. This gives the values of $\sgnvar\sigma$ near $\alpha$ and $\beta$. Since $|F| \to +\infty$ at $\alpha$ and $\beta$, the roots of $F$ lie in a compact subinterval of $I$, and Proposition \ref{prop:sigma} gives $n(F;I) \leq V^\gamma_\alpha(c)$ with an even difference. The last term $\sigma_s(\alpha)=c_s\delta_{s,s}p_s(\alpha)^{\gamma-s+1}$ has the sign of $c_s$, which gives the congruence $V^\gamma_\alpha(c) \equiv \sgnvar(c)$ as in the proof of Corollary \ref{cor:VleqD}, and $V^\gamma_\alpha(c) \leq \sgnvar(c)$ follows from the first item. Finally, since $\alpha=-a_1/b_1$ and $b_1>0$, we have $p_i(\alpha)=\frac{a_ib_1-a_1b_i}{b_1}=\frac{|D_{1i}|}{b_1}$, so that $\delta_{2,i}=|D_{1i}|=b_1p_i(\alpha)$ and $\sigma_2(\alpha)=b_1\sum_{i \geq 2}c_ip_i(\alpha)^\gamma$.

We prove the first item by induction on $s$, following the proof of Proposition \ref{prop:BFleqD}. If $s=1$ then $\sigma=(c_1p_1^\gamma)$ and $\sgnvar\sigma=0$. Assume $s \geq 2$ and fix $y \in I$.

\emph{Reduction to $I$ bounded.} After a translation we may assume $y=0$. With $\lambda$ and $\varphi(t)=\frac{t}{1-\lambda t}$ as in Step 1 of the proof of Proposition \ref{prop:BFleqD}, the linear functions $\tilde p_j$ of Lemma \ref{lem:projective} satisfy \eqref{hyp:L} on a bounded interval $\tilde I \ni 0$, and near $t=0$ we have $\tilde p_j(t)^\gamma=w(t)\,f_j(\varphi(t))$ with $w(t)=(1-\lambda t)^\gamma>0$. As in the proof of Lemma \ref{lem:projective}, it follows that $H_k(\tilde F)(t)=w(t)^k\varphi'(t)^{k(k-1)/2}H_k(F)(\varphi(t))$, where $\tilde F=\sum_jc_j\tilde p_j^\gamma$, hence $H_k(\tilde F)(0)=H_k(F)(0)$. By Proposition \ref{prop:sigma}, the signs of the entries of $\sigma(0)$ are thus the same for the two families, and we may assume that $I=(\alpha,\beta)$ is bounded.

\emph{Reduction to generic coefficients.} As in Step 2 of the proof of Proposition \ref{prop:BFleqD}, we may replace $c$ by a nearby vector $c'$ with the same signs pattern as in that proof, so that $\sgnvar(c')=\sgnvar(c)$ and $\sgnvar\sigma(0)$ does not decrease. Since $\sigma_k(\alpha)$ is a linear form in $c$ whose coefficient of $c_k$ is $\delta_{k,k}p_k(\alpha)^{\gamma-k+1}>0$, we may moreover assume that all $c_i$ and all $\sigma_k(\alpha)$, $k \geq 2$, are non-zero.

\emph{Comparison with $\alpha$.} By Proposition \ref{prop:sigma}, $\sgnvar\sigma$ is non-increasing on $I$, so that $\sgnvar\sigma(0) \leq V^\gamma_\alpha(c)$ by the computation of the limit at $\alpha$ above.

\emph{Induction.} The linear functions $p_2,\ldots,p_s$ satisfy \eqref{hyp:L} on $I_2=\{p_2>0,\ldots,p_s>0\} \ni \alpha$. Let $\sigma'=(\sigma'_1,\ldots,\sigma'_{s-1})$ be the sequence \eqref{eq:sigma} associated with $p_2,\ldots,p_s$ and $(c_2,\ldots,c_s)$. Using $|D_{1i}|=b_1p_i(\alpha)$, we get for $k=2,\ldots,s$
$$\sigma_k(\alpha)=\sum_{i=k}^s c_i\,|D_{1i}|\prod_{j=2}^{k-1}|D_{ji}|\;p_i(\alpha)^{\gamma-k+1}=b_1\sum_{i=k}^s c_i\prod_{j=2}^{k-1}|D_{ji}|\;p_i(\alpha)^{\gamma-(k-1)+1}=b_1\,\sigma'_{k-1}(\alpha).$$
Hence $V^\gamma_\alpha(c)=\sgnvar(c_1,\sigma'(\alpha))$, and by the induction hypothesis $\sgnvar\sigma'(\alpha) \leq \sgnvar(c_2,\ldots,c_s)$. All the numbers involved are non-zero and the last term $\sigma'_{s-1}(\alpha)$ has the sign of $c_s$. The same parity argument as in Step 4 of the proof of Proposition \ref{prop:BFleqD} gives $\sgnvar(c_1,\sigma'(\alpha)) \leq \sgnvar(c_1,\ldots,c_s)$, which concludes the proof.
\end{proof}

\begin{rem}
\label{rem:powers}
\begin{enumerate}
\item For $\gamma=-1$ we recover the functions $1/p_j$ of Example \ref{firstD}. As in Remark \ref{rem:comparisonDescartes}, the inequality $n(F;I) \leq V^\gamma_\alpha(c)$ refines Descartes' rule of signs: for $s=3$, we get $V^\gamma_\alpha(c)=\sgnvar\big(c_1,\ c_2p_2(\alpha)^\gamma+c_3p_3(\alpha)^\gamma,\ c_3\big)$, so that if $c_1,c_3>0>c_2$ and $c_3p_3(\alpha)^\gamma>|c_2|\,p_2(\alpha)^\gamma$, then $F$ has no root in $I$, while Descartes' rule of signs allows two roots.
\item The assumption $\gamma<0$ is essential for Theorem \ref{thm:powers}. For $\gamma \notin \{0,\ldots,s-2\}$, Lemma \ref{lem:powersW} shows that Theorem \ref{application} applies to $f_1,\ldots,f_s$ after multiplying them by suitable signs (with $f_1$ unchanged). However, the first item of Theorem \ref{thm:powers} fails in general for $\gamma>0$. For instance, take $\gamma=\frac{1}{2}$ and $p_1=y$, $p_2=1$, $p_3=1-y$ on $I=(0,1)$, which satisfy \eqref{hyp:L}, and $F=f_1+f_2+f_3$. Then $F>0$, while
$W(f_1,F)=W(f_1,f_2)+W(f_1,f_3)=-\frac{1}{2\sqrt{y}}-\frac{1}{2\sqrt{y(1-y)}}<0$ on $I$, so that the first two terms of the Budan--Fourier sign sequence already have one sign change, whereas $\sgnvar(1,1,1)=0$. Moreover, for $\gamma>0$ the terms $p_s^{\gamma-k+1}$ with $k-1<\gamma$ tend to $0$ at $\beta$, so that the sign sequence near $\beta$ is no longer determined by $c_s$ alone.
\item Many other families satisfying Descartes' rule of signs can be found in \cite[Part V]{P-S}, for instance exponential functions $e^{\lambda_1x},\ldots,e^{\lambda_sx}$ with $\lambda_1<\cdots<\lambda_s$ on $\R$ (Laguerre's rule of signs for exponential sums), monomials $x^{a_1},\ldots,x^{a_s}$ with real exponents on $(0,+\infty)$ (Example \ref{ex:monomialsW}), their confluent versions such as $x^je^{\lambda x}$ or $x^a(\log x)^j$, and more generally families $K(x,t_1),\ldots,K(x,t_s)$ with $t_1<\cdots<t_s$ defined by a totally positive kernel $K$ (see also \cite{Karlin}). By Proposition \ref{prop:renumbering}, Theorem \ref{application} applies to all of them. For exponentials and monomials, the functions have pairwise distinct orders of growth at both endpoints of the interval, and the limits of the Budan--Fourier sign sequence at the endpoints are $\sgnvar(c_1,\ldots,c_s)$ and $0$ (see Example \ref{ex:monomialsW}): near an endpoint, the terms of the sign sequence are dominated by the functions in the order of their growth, so that their signs are determined by the signs of the coefficients. Thus Theorem \ref{application} only localizes Descartes' rule of signs for these families. The situation is different in Corollary \ref{cor:VleqD} and Theorem \ref{thm:powers}: there, only the first function is singular at the left endpoint $\alpha$, while the other functions have finite values at $\alpha$. These values enter the limit of the sign sequence at $\alpha$, which explains why the resulting bound can be smaller than $\sgnvar(c_1,\ldots,c_s)$. It would be interesting to find other families in \cite[Part V]{P-S} with this property.
\end{enumerate}
\end{rem}

\section{Budan--Fourier for systems supported on circuits}
\label{sec:circuits}

\subsection{Recollection from \cite{BDF}}
\label{ssec:recollection}
Let $\A=\{a_0,\ldots,a_{n+1}\} \subset \Z^n$ be a circuit, that is, a minimally affinely dependent configuration whose convex hull has dimension $n$, and let $b=(b_0,\ldots,b_{n+1})$ be a non-zero affine relation on $\A$: $\sum_j b_j=0$ and $\sum_j b_ja_j=0$ (all $b_j$ are non-zero since $\A$ is a circuit). For a coefficient matrix $C=(c_{i,j}) \in \R^{n \times (n+2)}$ of rank $n$, consider the system
\begin{equation}\label{E:system}
f_i(x)=\sum_{j=0}^{n+1} c_{i,j}\,x^{a_j} = 0 \, , \quad i=1,\ldots,n,
\end{equation}
and denote by $n_\A(C)$ its number of positive solutions counted with multiplicities.

Let $P_0,\ldots,P_{n+1} \in \R^2$ be a Gale dual configuration of the columns of $C$, that is, the rows of a matrix whose columns form a basis of $\ker C$. Set $p_j(y)=\langle P_j,(1,y)\rangle$ and $\Delta_P=\{y \in \R \, | \, p_j(y)>0, \; j=0,\ldots,n+1\}$. Since $C$ has rank $n$, the matrix with rows $P_j$ has rank $2$, so the vectors $P_0,\ldots,P_{n+1}$ span $\R^2$.
If $(\R_{>0})^{n+2} \cap \ker C=\emptyset$ then $n_\A(C)=0$. Otherwise the open convex cone
$K=\{w \in \R^2 \, | \, \langle P_j,w\rangle>0, \; j=0,\ldots,n+1\}$ is non-empty, and the positive vectors of $\ker C$ are the vectors $(\langle P_j,w\rangle)_j$ with $w \in K$.

\begin{lemma}[{\cite[Lemma 2.3]{BDG}}]\label{lem:bounded}
Assume $K \neq \emptyset$. After a linear change of coordinates on $\R^2$ (that is, after replacing the Gale dual configuration $\{P_j\}$ by $\{P_jM\}$ for some $M \in GL_2(\R)$), the first coordinate is positive on $\overline{K}\setminus\{0\}$. Then every $w \in K$ is a positive multiple of a unique vector $(1,y)$ with $y \in \Delta_P$, and $\Delta_P$ is a non-empty bounded open interval.
\end{lemma}

\begin{proof}
We include a short proof for the reader's convenience. The linear form $\ell(w)=\langle P_0+\cdots+P_{n+1},w\rangle$ is non-negative on $\overline{K}$, and if $\ell(w)=0$ for some $w \in \overline{K}$ then $\langle P_j,w\rangle=0$ for all $j$, hence $w=0$ since the $P_j$ span $\R^2$. Thus $\ell>0$ on $\overline{K}\setminus\{0\}$, and we choose coordinates in which $\ell$ is the first coordinate. Then every $w \in K$ is a positive multiple of a unique $(1,y)$, and $y \in \Delta_P$. Finally, $\Delta_P$ is the intersection of $K$ with the affine line $\{w_1=1\}$, which is bounded since $\overline{K}\setminus\{0\}$ is contained in $\{w_1>0\}$ and $\overline{K}\cap\{w_1=1\}$ is compact.
\end{proof}

In the rest of the paper we always choose the Gale dual configuration as in Lemma \ref{lem:bounded}, so that $\Delta_P$ is a non-empty bounded open interval.

\begin{prop}[{\cite{BS07}, see also \cite[\S 4]{BD}}]
\label{prop:Galebijection}
Let $g(y)=\prod_{j} p_j(y)^{b_j}$. The map sending a positive solution $x$ of \eqref{E:system} to the unique $y \in \Delta_P$ such that $(x^{a_0},\ldots,x^{a_{n+1}})$ is proportional to $(p_0(y),\ldots,p_{n+1}(y))$ is a bijection onto the set of solutions of $g(y)=1$ in $\Delta_P$, which preserves multiplicities.
\end{prop}

Group the indices $j$ according to the rays $\R_{>0}P_j$, and order these $k$ rays by a strict ordering $\bar\sigma\colon\{0,\ldots,k-1\} \to \{0,\ldots,n+1\}$ (choosing a representative $P_{\bar\sigma_\ell}$ on each ray) so that $\det(P_{\bar\sigma_i},P_{\bar\sigma_j})$ has a constant sign for $i<j$, see \cite{BDF}. Let $\lambda_\ell$ be the sum of the $b_j$ over the indices $j$ on the $\ell$-th ray, and $\mu_\ell=\lambda_0+\cdots+\lambda_\ell$, so that $\mu_{k-1}=0$. Since $p_j$ is a positive multiple of $p_{\bar\sigma_\ell}$ when $P_j$ lies on the $\ell$-th ray, there is a constant $e>0$ such that
\begin{equation}\label{eq:gbis}
g = e \cdot \prod_{\ell=0}^{k-2} \left( \frac{p_{\bar{\sigma}_\ell}}{p_{\bar{\sigma}_{\ell+1}}}\right)^{\mu_\ell}, \qquad
G'=(\log g)' = \sum_{\ell=0}^{k-2} \mu_\ell
\left(\frac{p'_{\bar{\sigma}_\ell}}{p_{\bar{\sigma}_\ell}} - \frac{p'_{\bar{\sigma}_{\ell+1}}}{p_{\bar{\sigma}_{\ell+1}}}\right) \quad \text{on } \Delta_P.
\end{equation}
The main results of \cite{BDF} are the following.

\begin{thm}[{\cite{BDF}}]\label{thm:refines}
If $n_\A(C)$ is finite, then $n_\A(C) \leq 1+\sgnvar(\mu)$. Moreover, this bound is optimal: for any circuit $\A$ there exists $C$ such that equality holds.
\end{thm}

\begin{prop}[{\cite{BDF}}]\label{prop:congruence}
If moreover $\lambda_0 \neq 0$ and $\lambda_{k-1} \neq 0$, then $n_\A(C) \equiv 1+\sgnvar(\mu) \mod 2$.
\end{prop}

\subsection{Comparison with the bound of \cite{BDF}}

\begin{rem}
\label{rem:comparisonBDF}
In \cite{BDF}, Theorem \ref{thm:refines} is proved as follows. By \eqref{eq:gbis}, the derivative of $G=\log g$ is the linear combination
$G'=\sum_{\ell=0}^{k-2}\mu_\ell\, h_{\ell}$ of the functions $h_\ell=\frac{p'_{\bar\sigma_\ell}}{p_{\bar\sigma_\ell}}-\frac{p'_{\bar\sigma_{\ell+1}}}{p_{\bar\sigma_{\ell+1}}}$, which are of the type considered in Example \ref{secD} (with $q_j=\frac{p'_j}{p_j}$), provided that $\det(P_{\bar\sigma_i},P_{\bar\sigma_j})<0$ for $i<j$. This can always be achieved by reversing the strict ordering $\bar\sigma$ if necessary, which replaces $(\mu_0,\ldots,\mu_{k-2})$ by $(-\mu_{k-2},\ldots,-\mu_0)$ and does not change $\sgnvar(\mu)$. These functions satisfy Descartes' rule of signs on $\Delta_P$, so that $G'$ has at most $\sgnvar(\mu)$ roots in $\Delta_P$, and Rolle's theorem gives $n_\A(C) \leq 1+\sgnvar(\mu)$.
\end{rem}

We now apply Theorem \ref{application} directly to the function $G=\log g$, instead of applying Descartes' rule of signs to $G'$ and then Rolle's theorem. This gives bounds for the number of positive solutions of the circuit system \eqref{E:system} whose logarithms lie in some explicit polyhedra, which refine the bound of \cite{BDF} and contain the parity statement of Proposition \ref{prop:congruence}.

\smallskip

\noindent\emph{Assumptions and notations.} We use the notations of Section \ref{ssec:recollection}, and we assume:
\begin{equation}
\label{hyp:C}
\tag{$\mathcal{C}$}
\begin{minipage}{0.85\textwidth}
$\A$ is a circuit, $C$ has rank $n$, $K \neq \emptyset$ and the Gale dual configuration is chosen as in Lemma \ref{lem:bounded} (so $\Delta_P=(\alpha,\beta)$ is non-empty and bounded), the strict ordering $\bar\sigma$ satisfies $\det(P_{\bar\sigma_i},P_{\bar\sigma_j})<0$ for $i<j$, and $\mu \neq 0$.
\end{minipage}
\end{equation}
The condition on $\bar\sigma$ can always be achieved by reversing the strict ordering (Remark \ref{rem:comparisonBDF}), and if $\mu=0$ then $g=e$ is constant, so that $n_\A(C)=0$ or $n_\A(C)=\infty$. We set:
\begin{itemize}
\item $N=1+\max\{\ell \, | \, \mu_\ell \neq 0\}$, and for $i=1,\ldots,N$,
$c_i=\mu_{i-1}$, $\psi_i=\log\frac{p_{\bar\sigma_{i-1}}}{p_{\bar\sigma_{i}}}$ and
$h_i=\psi_i'=\frac{p'_{\bar\sigma_{i-1}}}{p_{\bar\sigma_{i-1}}}-\frac{p'_{\bar\sigma_{i}}}{p_{\bar\sigma_{i}}}$, so that by \eqref{eq:gbis}
$$G=\log g=\log e+\sum_{i=1}^N c_i\psi_i, \qquad G'=\sum_{i=1}^N c_ih_i, \qquad c_N \neq 0;$$
\item $H_m(G')=W(h_1,\ldots,h_{m-1},G')$ for $m=1,\ldots,N$, and for $y \in \Delta_P$
\begin{equation}
\label{eq:S}
S(y)=\sgnvar\big(G(y),\, H_1(G')(y),\, H_2(G')(y),\ldots,H_N(G')(y)\big);
\end{equation}
\item for an interval $J \subset \Delta_P$, $n_\A(C;J)$ is the number of positive solutions $x$ of \eqref{E:system} whose associated point $y(x) \in \Delta_P$ (Proposition \ref{prop:Galebijection}) lies in $J$, counted with multiplicities.
\end{itemize}
Note that the linear functions $p_{\bar\sigma_0},\ldots,p_{\bar\sigma_N}$ satisfy \eqref{hyp:L} (with $s=N+1$) on an interval containing $\Delta_P$, and that $h_1,\ldots,h_N$ are the corresponding functions of Section \ref{ssec:setting3}.

The number $n_\A(C;J)$ can be described without referring to $y$, in logarithmic coordinates. For $x \in \R^n_{>0}$ write $\log x=(\log x_1,\ldots,\log x_n)$. For $i,j \in \{0,\ldots,n+1\}$ let $\rho_{ij}=\frac{p_i}{p_j}$. If $P_i$ and $P_j$ are collinear, then $\rho_{ij}$ is a positive constant on $\Delta_P$; otherwise $\rho_{ij}'=\frac{\det(P_j,P_i)}{p_j^2}$ does not vanish, and $\rho_{ij}$ is strictly monotone on $\Delta_P$. For an interval $J=(a,b] \subset \Delta_P$, the set $\rho_{ij}(J)$ is thus either a point or a half-open interval with endpoints $\rho_{ij}(a)$ (excluded) and $\rho_{ij}(b)$ (included). We set
\begin{equation}
\label{eq:PJ}
\mathcal{P}_J=\Big\{z \in \R^n \; \Big| \; \langle a_i-a_j,z\rangle \in \log \rho_{ij}(J) \text{ for all } i,j \in \{0,\ldots,n+1\}\Big\}.
\end{equation}
Each condition in \eqref{eq:PJ} defines either a set of the form $\{z \in \R^n \, | \, c_1<\langle a_i-a_j,z\rangle \leq c_2\}$ (a half-open slab between two parallel affine hyperplanes), when $P_i$ and $P_j$ are not collinear, or an affine hyperplane $\{\langle a_i-a_j,z\rangle=c\}$, when they are collinear. Thus $\mathcal{P}_J$ is a (half-open) polyhedron. If moreover $\alpha<a<b<\beta$, then each $\rho_{i0}$ is bounded and bounded away from $0$ on $J$, and since the vectors $a_i-a_0$, $i=1,\ldots,n+1$, span $\R^n$ (the convex hull of $\A$ has dimension $n$), the polyhedron $\mathcal{P}_J$ is bounded, that is, a (half-open) polytope.

\begin{lemma}
\label{lem:PJ}
Let $J=(a,b] \subset \Delta_P$. A positive solution $x$ of \eqref{E:system} satisfies $y(x) \in J$ if and only if $\log x \in \mathcal{P}_J$. Consequently, $n_\A(C;J)$ is the number of positive solutions $x$ of \eqref{E:system} such that $\log x \in \mathcal{P}_J$, counted with multiplicities. Moreover, the polyhedra $\mathcal{P}_{(y_{t-1},y_t]}$ associated with a subdivision $y_0<y_1<\cdots<y_r$ of an interval contained in $\Delta_P$ are pairwise disjoint.
\end{lemma}

\begin{proof}
Since $(x^{a_0},\ldots,x^{a_{n+1}})$ is proportional to $(p_0(y(x)),\ldots,p_{n+1}(y(x)))$, we have $\langle a_i-a_j,\log x\rangle=\log x^{a_i-a_j}=\log\rho_{ij}(y(x))$ for all $i,j$. Hence $y(x) \in J$ implies $\log x \in \mathcal{P}_J$. Conversely, since the vectors $P_j$ span $\R^2$, there exist $i,j$ such that $P_i$ and $P_j$ are not collinear; then $\rho_{ij}$ is injective on $\Delta_P$, and $\log\rho_{ij}(y(x)) \in \log\rho_{ij}(J)$ implies $y(x) \in J$. The last assertion follows since, for such a pair $(i,j)$, the sets $\rho_{ij}((y_{t-1},y_t])$ are pairwise disjoint.
\end{proof}

\begin{thm}
\label{thm:localized}
Assume \eqref{hyp:C} and let $S$ be defined by \eqref{eq:S}.
\begin{enumerate}
\item For any $a,b \in \Delta_P$ with $a<b$, we have
$$n_\A(C;(a,b]) \leq S(a)-S(b),$$
and the difference $S(a)-S(b)-n_\A(C;(a,b])$ is an even integer. By Lemma \ref{lem:PJ}, $n_\A(C;(a,b])$ is the number of positive solutions $x$ of \eqref{E:system} with $\log x$ in the polyhedron $\mathcal{P}_{(a,b]}$ defined in \eqref{eq:PJ}.
\item For any $y \in \Delta_P$ we have $S(y) \leq 1+\sgnvar(\mu)$. In particular, if $\alpha=y_0<y_1<\cdots<y_r=\beta$ is a subdivision of $\Delta_P=(\alpha,\beta)$, then the bounds of item (1) for the intervals $(y_{t-1},y_t]$, $t=2,\ldots,r-1$, which partition $(y_1,y_{r-1}]$, add up to $S(y_1)-S(y_{r-1}) \leq 1+\sgnvar(\mu)$.
\item Assume moreover that $\lambda_0 \neq 0$ and $\lambda_{k-1} \neq 0$ (that is, $\mu_0 \neq 0$ and $\mu_{k-2} \neq 0$). Then $S(y)=0$ for $y$ close enough to $\beta$, the limit $S(\alpha^+)=\lim_{y \to \alpha^+}S(y)$ exists, and
$$
n_\A(C) \leq S(\alpha^+) \leq 1+\sgnvar(\mu), \qquad n_\A(C) \equiv S(\alpha^+) \equiv 1+\sgnvar(\mu) \mod 2.
$$
If moreover $T_k(c) \neq 0$ for $k=1,\ldots,N-1$ (see \eqref{eq:Valpha}), then $S(\alpha^+)=1+V_\alpha(c)$.
\end{enumerate}
\end{thm}

\begin{proof}
The linear functions $p_{\bar\sigma_0},\ldots,p_{\bar\sigma_N}$ satisfy the hypotheses of Example \ref{secD} on the interval $I'=\{y \in \R \, | \, p_{\bar\sigma_0}(y)>0,\ldots,p_{\bar\sigma_N}(y)>0\}$, which contains $\Delta_P$ (with $s=N+1$ and $q_{i}=\frac{p'_{\bar\sigma_{i-1}}}{p_{\bar\sigma_{i-1}}}$, $i=1,\ldots,N+1$). Thus $W(h_1,\ldots,h_m)>0$ on $\Delta_P$ for $m=1,\ldots,N$.

Consider the family of $N+1$ analytic functions $(\varphi_1,\ldots,\varphi_{N+1})=(1,\psi_1,\ldots,\psi_N)$ on $\Delta_P$. Expanding the Wronskian along its first column $(1,0,\ldots,0)^T$, we get for any analytic functions $f_1,\ldots,f_m$ the identity
\begin{equation}
\label{eq:W1}
W(1,f_1,\ldots,f_m)=W(f_1',\ldots,f_m').
\end{equation}
Hence $W(\varphi_1)=1$ and $W(\varphi_1,\ldots,\varphi_{m+1})=W(h_1,\ldots,h_m)>0$ on $\Delta_P$ for $m=1,\ldots,N$, so that Theorem \ref{application} applies to the family $(\varphi_1,\ldots,\varphi_{N+1})$ on $\Delta_P$ and to $G=\log e\cdot\varphi_1+\sum_{i=1}^N c_i\varphi_{i+1}$, whose last coefficient $c_N$ is non-zero. By \eqref{eq:W1}, the corresponding sequence is
$$
H_1(G)=G, \qquad H_{m+1}(G)=W(1,\psi_1,\ldots,\psi_{m-1},G)=W(h_1,\ldots,h_{m-1},G')=H_m(G'), \quad m=1,\ldots,N,
$$
so that $\sgnvar H(G)(y)=S(y)$. Theorem \ref{application} gives that the number of roots of $G$ in $(a,b]$, counted with multiplicities, is at most $S(a)-S(b)$, with an even difference. Since $G=\log g$, the roots of $G$ in $\Delta_P$ are the solutions of $g=1$ in $\Delta_P$ with the same multiplicities, and the correspondence of Proposition \ref{prop:Galebijection} between positive solutions $x$ of \eqref{E:system} and solutions $y$ of $g(y)=1$ in $\Delta_P$ is local in $y$. This proves item (1).

For item (2), the sequence defining $S(y)$ is $(G(y),H(G')(y))$, so $S(y) \leq 1+\sgnvar H(G')(y)$. By Proposition \ref{prop:BFleqD} applied to $h_1,\ldots,h_N$ on $I'$, we have $\sgnvar H(G')(y) \leq \sgnvar(c_1,\ldots,c_N)=\sgnvar(\mu)$. The second assertion of item (2) follows since the bounds telescope.

Let us prove item (3). Now $N=k-1$, the interval $I'$ is equal to $\Delta_P$, and $\alpha$ (respectively $\beta$) is the root of $p_{\bar\sigma_0}$ (respectively $p_{\bar\sigma_{N}}$), see the discussion before Proposition \ref{prop:endpoints}.
When $y \to \beta^-$, we have $G(y)=-c_N\log p_{\bar\sigma_N}(y)+O(1)$, so that $G$ has the sign of $c_N$ near $\beta$; by item (1) of Proposition \ref{prop:endpoints} (whose proof only uses $c_N \neq 0$), all $H_m(G')$ have the sign of $c_N$ near $\beta$. Thus $S(y)=0$ near $\beta$.
By Theorem \ref{application}, $S$ is non-increasing on $\Delta_P$ and integer-valued, so $S(\alpha^+)$ exists.
Since $G(y) \to \pm \infty$ when $y$ tends to $\alpha$ or $\beta$, the roots of $G$ lie in a compact subset of $\Delta_P$. Taking $a$ close to $\alpha$ and $b$ close to $\beta$ in item (1), we get $n_\A(C) \leq S(\alpha^+)$ with an even difference, and $S(\alpha^+) \leq 1+\sgnvar(\mu)$ by item (2).
When $y \to \alpha^+$, we have $G(y)=c_1\log p_{\bar\sigma_0}(y)+O(1)$, so that $G$ has the sign of $-c_1=-\mu_0$ near $\alpha$, while the last term $H_N(G')=c_NW(h_1,\ldots,h_N)$ of the sequence has the sign of $c_N=\mu_{k-2}$. Since the number of sign changes of a sequence of real numbers whose first and last terms are non-zero has the parity of the number of sign changes between these two terms, we get $S(\alpha^+) \equiv \sgnvar(-\mu_0,\mu_{k-2}) \equiv 1+\sgnvar(\mu_0,\mu_{k-2}) \equiv 1+\sgnvar(\mu) \mod 2$.
Finally, near $\alpha$ we have $H_1(G')=G'=c_1q_1+O(1)$ with $q_1 \to +\infty$, so $G$ and $H_1(G')$ have opposite signs, and the remaining terms are given by item (2) of Proposition \ref{prop:endpoints} under the genericity assumption. This gives $S(\alpha^+)=1+V_\alpha(c)$.
\end{proof}

\begin{cor}
\label{cor:compact}
Assume \eqref{hyp:C} and let $S$ be defined by \eqref{eq:S}. Let $i_0=\bar\sigma_0$ and $i_1=\bar\sigma_{k-1}$, and consider the monomial $m(x)=x^{a_{i_0}-a_{i_1}}$ and the function $\rho=\frac{p_{i_0}}{p_{i_1}}$, which is an increasing bijection from $\Delta_P$ onto $(0,+\infty)$. Let $K \subset \R^n_{>0}$ be a compact set, let $m_-=\min_K m$, $m_+=\max_K m$, and $a'=\rho^{-1}(m_-)$, $b'=\rho^{-1}(m_+)$, so that $\alpha<a' \leq b'<\beta$. Then the number of positive solutions of \eqref{E:system} in $K$, counted with multiplicities, is at most
$$S(a'^-)-S(b') \leq 1+\sgnvar(\mu),$$
where $S(a'^-)=\lim_{y \to a'^-}S(y)$. Moreover, the number of positive solutions of \eqref{E:system} in the region $R_K=\{x \in \R^n_{>0} \, | \, m_- \leq m(x) \leq m_+\}$, which contains $K$, is congruent to $S(a'^-)-S(b')$ modulo $2$, and all these solutions $x$ satisfy $\log x \in \mathcal{P}_{(a,b']}$ for any $a \in (\alpha,a')$; in particular they lie in a bounded subset of $\R^n_{>0}$.
\end{cor}

\begin{proof}
By \eqref{hyp:C}, $\det(P_{i_0},P_{i_1})<0$, so that $\rho'=\frac{\det(P_{i_1},P_{i_0})}{p_{i_1}^2}>0$ on $\Delta_P$. The linear functions $p_{\bar\sigma_0},\ldots,p_{\bar\sigma_{k-1}}$ satisfy \eqref{hyp:L} on $\Delta_P=(\alpha,\beta)$, so that $p_{i_0}(\alpha)=0$ and $p_{i_1}(\beta)=0$ (see Subsection \ref{ssec:setting3}); hence $\rho$ tends to $0$ at $\alpha$ and to $+\infty$ at $\beta$, and $\rho$ is an increasing bijection from $\Delta_P$ onto $(0,+\infty)$. Since $K$ is compact, $0<m_- \leq m_+<+\infty$, and thus $\alpha<a' \leq b'<\beta$.

For a positive solution $x$, we have $m(x)=\rho(y(x))$ since $(x^{a_0},\ldots,x^{a_{n+1}})$ is proportional to $(p_0(y(x)),\ldots,p_{n+1}(y(x)))$. Hence $x \in R_K$ if and only if $y(x) \in [a',b']$. By Theorem \ref{thm:localized}, $S$ is non-increasing on $\Delta_P$, so the limit $S(a'^-)$ exists, and by item (1) of Theorem \ref{thm:localized} there are finitely many solutions $y$ of $g(y)=1$ in any interval $(a,b'] $ with $a>\alpha$. Choose $a \in (\alpha,a')$ such that $S(a)=S(a'^-)$ and no solution $y(x)$ lies in $(a,a')$. Then the positive solutions in $R_K$ are exactly those with $y(x) \in (a,b']$, and item (1) of Theorem \ref{thm:localized} shows that their number is at most $S(a)-S(b')=S(a'^-)-S(b')$, with an even difference. Since $K \subset R_K$, this bounds the number of positive solutions in $K$. By item (2) of Theorem \ref{thm:localized}, $S(a'^-) \leq 1+\sgnvar(\mu)$, and $S(b') \geq 0$. Finally, by Lemma \ref{lem:PJ}, every positive solution $x$ with $y(x) \in (a,b']$ satisfies $\log x \in \mathcal{P}_{(a,b']}$, which is bounded since $\alpha<a<b'<\beta$.
\end{proof}

\begin{rem}
\label{rem:localized}
\begin{enumerate}
\item The bound of \cite{BDF} depends only on $\A$ and on the strict ordering $\bar\sigma$, that is, on the signs of the maximal minors of $C$. The bounds of Theorem \ref{thm:localized} depend in addition on the values of $G$ and of the Wronskians $H_m(G')$ at some points of $\Delta_P$, which are expressed in terms of the Gale dual vectors $P_j$, that is, of ratios of maximal minors of $C$. Thus Theorem \ref{thm:localized} is not a rule of signs, but a refinement of the bound of \cite{BDF} inside each set of coefficient matrices with a given strict ordering.
\item This does not contradict the optimality of the bound of \cite{BDF} (Theorem \ref{thm:refines}): under the assumptions of item (3) of Theorem \ref{thm:localized}, a matrix $C$ such that $n_\A(C)=1+\sgnvar(\mu)$ necessarily satisfies $S(\alpha^+)=1+\sgnvar(\mu)$.
\item For $k \leq 3$ we have $N \leq 2$, and $S(\alpha^+)=1+\sgnvar(\mu)$ under the assumptions of item (3). For $k=4$ (for instance, for uniform matrices $C$ when $n=2$) and $\mu_0\mu_2 \neq 0$, we have $N=3$, $c=(\mu_0,\mu_1,\mu_2)$ and, by Remark \ref{rem:comparisonDescartes},
$V_\alpha(c)=\sgnvar\big(\mu_0,\ \mu_1h_2(\alpha)+\mu_2h_3(\alpha),\ \mu_2\big)$. Hence, if $\sgnvar(\mu)=2$ (so that \cite{BDF} allows three positive solutions) and $\mu_1h_2(\alpha)+\mu_2h_3(\alpha)$ has the sign of $\mu_0$, then $n_\A(C) \leq 1$. This describes an explicit semi-algebraic set of coefficient matrices with at most one positive solution, which is smaller than the sets given by the strict orderings.
\end{enumerate}
\end{rem}

\section*{Declaration of generative AI and AI-assisted technologies in the writing process}
Claude Opus 5.5 (Anthropic) helped the author in the comparison with Descartes' rule of signs (Propositions \ref{prop:endpoints} and \ref{prop:BFleqD} and Corollary \ref{cor:VleqD}) and in the editing of the text. It also helped in obtaining Theorem \ref{thm:powers} on negative powers of linear functions. The author takes full responsibility for the contents of this paper.

\providecommand{\bysame}{\leavevmode\hbox to3em{\hrulefill}\thinspace}
\providecommand{\MR}{\relax\ifhmode\unskip\space\fi MR }
\providecommand{\MRhref}[2]{%
  \href{http://www.ams.org/mathscinet-getitem?mr=#1}{#2}
}
\providecommand{\href}[2]{#2}

\end{document}